\documentclass[11pt,a4paper,reqno]{amsart}
\usepackage[margin=1.1in]{geometry}
\usepackage{amssymb,amsmath,graphicx,amsfonts}
\usepackage{mathrsfs}
\usepackage{multirow}
\usepackage{graphicx,epsfig}
\usepackage{color}
\usepackage{enumerate,enumitem}
\usepackage{empheq}
\usepackage{cases}
\usepackage{cite}
\usepackage{mathrsfs}
\usepackage{amsgen}
\usepackage{amscd}

\allowdisplaybreaks

\usepackage{booktabs}
\usepackage{amsmath}
\usepackage{latexsym}
\usepackage{amssymb,amsmath,graphicx,amsfonts,euscript}
\usepackage{amsfonts,bm,hyperref,subcaption}
\usepackage{enumerate}
\newtheorem{theorem}{Theorem}[section]

\newtheorem{remark}[theorem]{Remark}
\theoremstyle{definition}

\newtheorem{example}[theorem]{Example}

\numberwithin{theorem}{section}
\numberwithin{equation}{section}

\def\le{\leqslant}
\def\ge{\geqslant}
\def\Omega{\varOmega}
\def\Delta{\varDelta}
\def\bex{\begin{exercise}\upshape}
\def\eex{\end{exercise}}

\numberwithin{equation}{section}
\begin{document}

\title[]{An energy-stable minimal deformation rate scheme for mean curvature flow and surface diffusion}\thanks{This work is supported in part by AMSS-PolyU Joint Laboratory and the Research Grants Council of Hong Kong (PolyU/GRF15304024 and PolyU/RFS2324-5S03).}

\author[]{Guangwei Gao,\,\, Harald Garcke,\,\, Buyang Li,\,\, and\,\,Rong Tang}
\address{Guangwei Gao, Buyang Li, Rong Tang: Department of Applied Mathematics, The Hong Kong Polytechnic University, Hong Kong. 
{\rm Email address: {\tt guang-wei.gao@polyu.edu.hk}, {\tt buyang.li@polyu.edu.hk} and {\tt claire.tang@polyu.edu.hk}
}}

\address{Harald Garcke: Fakult\"at f\"ur Mathematik, Universität Regensburg, 93040 Regensburg, Germany. {\rm Email address: {\tt harald.garcke@ur.de}
}}




\maketitle

\begin{abstract} 
We propose a new parametric finite element method, referred to as the BGN-MDR method, for simulating both mean curvature flow and surface diffusion for closed hypersurfaces, as well as open hypersurfaces with moving contact lines in three dimensions. The method is also applicable to closed and open curves with moving contact points in two dimensions. The proposed scheme inherits the energy stability from the BGN scheme proposed by Barrett, Garcke, and N\"urnberg in 2008, and offers improved mesh quality similar to the minimal deformation rate (MDR) method proposed by Hu and Li in 2022, especially for small time step sizes where the BGN scheme may become unstable and result in deteriorated meshes. 
\end{abstract}


\setlength\abovedisplayskip{4pt}
\setlength\belowdisplayskip{4pt}

\section{Introduction}
The evolution of surfaces under geometric curvature flows, such as mean curvature flow and surface diffusion, has attracted significant interest due to its wide-ranging applications (see, e.g., \cite{BarrettGarckeNurnberg2020,deckelnick2005computation,
	ecker2012regularity,
	ganesan2017ale,li2021convergence}). These curvature-driven mechanisms are fundamental to understanding  interface dynamics in physical systems, as illustrated by phenomena such as grain boundary migration \cite{Gottstein2009} and the morphological evolution of thin solid films \cite{Thompson1}. 

This article concerns the numerical approximation of a  hypersurface \(\Gamma(t)\), \(t \in [0, T]\), in $\mathbb{R}^{d}$, where $d=2$ or $3$, evolving under geometric curvature flows such as mean curvature flow and surface diffusion, which are described by the following geometric evolution equations:
\begin{align*}
\begin{aligned}
v \cdot n &= -H && \text{(mean curvature flow)}, \\
v \cdot n &= \Delta_{\Gamma(t)} H && \text{(surface diffusion)}.
\end{aligned}
\end{align*}
Here, \(v(\cdot,t)\) denotes the velocity field of \(\Gamma(t)\), and the specific velocity law—whether corresponding to mean curvature flow or surface diffusion—governs the evolution of the surface.

The numerical approximation of geometric curvature flows governed by various velocity laws has been extensively studied by the parametric finite element method (FEM), which was introduced by Dziuk \cite{Dziuk} in 1990 and has since been widely used for computing a variety of curvature flows; see, e.g., \cite{bonito2010parametric,dziuk2008computational,banesch2005finite}. However, computing mean curvature flow or surface diffusion by the parametric FEM poses significant challenges, particularly in preserving the mesh quality as the surface undergoes large deformations; nodal points tend to cluster and the mesh becomes increasingly distorted, potentially leading to computational breakdowns. Consequently, to address these issues, advanced mesh remeshing techniques have been developed (see, e.g.,~\cite{marchandise2011highquality, remacle2010highquality}). These methods dynamically reallocate mesh points and restore mesh quality whenever it falls below a prescribed threshold, thereby ensuring that the numerical computations are consistently performed on meshes with good quality.

An alternative to the remeshing techniques is to impose an artificial tangential motion that continuously improves the mesh quality of the evolving surface. This idea was introduced by Barrett, Garcke, and N\"urnberg in \cite{barrett2007parametric, barrett2008hypersurfaces, barrett2008willmore}, in which they designed a weak formulation that enforces the map \(X_h^m:\Gamma_h^{m-1}\rightarrow\Gamma_h^m\) between consecutive approximate surfaces to be discrete harmonic. As a result, the method helps to preserve the shape of mesh triangles to a certain degree. Specifically, the Barrett--Garcke--N\"urnberg (BGN) method for mean curvature flow can be formulated as follows: Find \((v_h^m, H_h^m)\in (S_h^{m-1})^d \times S_h^{m-1}\), with $S_h^{m-1}$ being the piecewise linear finite element space on the surface $\Gamma_h^{m-1}$, such that
\begin{subequations}\label{BGN-weak-eqn}
    \begin{align}
        \int_{\Gamma_h^{m-1}}^{(h)} v_h^m \cdot n_h^{m-1}\,\chi_h &= \int_{\Gamma_h^{m-1}}^{(h)} -H_h^m\,\chi_h &&\forall\, \chi_h \in S_h^{m-1}, \label{BGN-weak-eqn_a} \\
        \int_{\Gamma_h^{m-1}} \nabla_{\Gamma_h^{m-1}} \Bigl(\tau v_h^m + \operatorname{id}\Bigr) \cdot \nabla_{\Gamma_h^{m-1}} \eta_h &= \int_{\Gamma_h^{m-1}}^{(h)} H_h^m\, n_h^{m-1} \cdot \eta_h && \forall\, \eta_h \in (S_h^{m-1})^d, \label{BGN-weak-eqn_b}
    \end{align}
\end{subequations}
where~\eqref{BGN-weak-eqn_b} enforces that \(X_h^m := \operatorname{id} + \tau v_h^m\) is a discrete harmonic map from \(\Gamma_h^{m-1}\) to \(\Gamma_h^m\), and \(\int_{\Gamma_h^{m-1}}^{(h)}\) denotes the mass-lumped quadrature for integration over \(\Gamma_h^{m-1}\). In addition, to enhancing mesh quality, the BGN method is also energy stable; that is, the surface area decreases over time. These two advantages (improving mesh quality while being energy stable) make the BGN method particularly suitable for computing mean curvature flow undergoing large deformations, often up to the onset of singularities—offering advantages over Dziuk's original parametric FEM. Beyond mean curvature flow, the BGN method has been successfully extended to a variety of applications, including two-phase Navier–Stokes flows~\cite{barrett2013eliminating, barrett2015stable}, surface diffusion~\cite{bao2021structure}, and axisymmetric geometric evolution equations~\cite{bao2022volume}. 

One drawback of the BGN method lies in the tangential velocity introduced in~\eqref{BGN-weak-eqn_b}. As \(\tau \to 0\), the tangential component becomes indeterminate, potentially leading to an ill-posed system. This issue may cause numerical instabilities or even breakdown of the simulation, particularly in three-dimensional settings.

In contrast to the BGN method, Elliott and Fritz~\cite{elliott2017approximations, elliott2016algorithms} introduced an artificial tangential velocity derived from reparametrizing the evolving surface using DeTurck flow techniques, which naturally yields the tangential component of the velocity field. A key advantage of this approach is that it facilitates rigorous convergence analysis of certain parametric FEMs incorporating tangential motion. This has been demonstrated through established convergence results for curve shortening flow~\cite{elliott2017approximations} and the mean curvature flow of closed torus-type surfaces~\cite{mierswa2020error}. However, the convergence of such methods for the mean curvature flow of general surfaces in three dimensions remains an open and compelling question. 


To address the instability of the BGN method in the limit \(\tau \to 0\), Hu and Li~\cite{hu2022evolving} (see also~\cite{bai2024convergent}) proposed an artificial tangential motion designed to minimize the energy 
\begin{align}\label{DR-energy}
\int_{\Gamma} |\nabla_\Gamma v|^2 
\end{align}
subject to the constraint \(v \cdot n = u \cdot n\), where \(u\) denotes the original velocity of the surface (for example $u=-Hn$ and $(\Delta_\Gamma H)n$ in mean curvature flow and surface diffusion, respectively). 
The energy in~\eqref{DR-energy} represents the instantaneous deformation rate of the surface. Minimizing this energy under the constraint \(v \cdot n = u \cdot n\) leads to the following continuous-level problem that determines the tangential component of the velocity: 
\begin{align}\label{MDR-eqn}
\begin{aligned}
-\Delta_{\Gamma} v &= \kappa\, n, \\
v \cdot n &= u \cdot n,
\end{aligned}
\end{align}
where \(\kappa\) appears as the Lagrange multiplier associated with the constrained minimization problem. 
{This approach, referred to as the \emph{minimal deformation rate} (MDR) method, demonstrates a significant distinction depending on the discretization used. Specifically, the MDR and BGN methods are equivalent under semi-discretization in time (i.e., without spatial discretization), as formally shown in \cite{hu2022evolving}. However, when spatial discretization is introduced, the two methods are no longer equivalent. In the limit $\tau \to 0$, corresponding to semi-discretization in space, the MDR method discretizes \eqref{MDR-eqn}, which uniquely determines the tangential motion. In contrast, the BGN method discretizes $\Delta_{\Gamma} {\rm id} = -H n$, an equation that holds for any surface and thus allows for arbitrary tangential motion. This fundamental difference explains the improved stability properties of the MDR method as $\tau \to 0$ in the fully discrete setting.
}
Despite its effectiveness in maintaining the mesh quality, a major challenge of the MDR method lies in constructing a linearly implicit full discretization scheme that preserves  energy stability—i.e., the area-decreasing property over time—that is inherent to the BGN method. Addressing this issue remains an important open problem in the development of the MDR approach. 

In addition to the MDR method, Duan and Li~\cite{duan2024new} proposed an alternative strategy that incorporates artificial tangential velocity by minimizing a deformation energy defined on the initial surface \(\Gamma_h^0\), while preserving the prescribed normal velocity. The continuous formulation of this \emph{minimal deformation} (MD) approach ensures that the flow map from the initial surface to the current surface is a harmonic map. As a result, this method can reduce mesh distortion and mitigate the accumulation of errors across successive time steps. Similar to the MDR approach, a major challenge in the MD approach lies in designing a linearly implicit, fully discrete scheme that preserves the energy stability—a key property of the BGN method. While this can be achieved through nonlinearly implicit schemes involving Lagrange multipliers~\cite{Gao-Li-2025}, constructing an efficient linearly implicit energy stable scheme that incorporates the MD tangential motion remains an open and interesting problem. 

%

Apart from advances in algorithm design, the convergence analysis of parametric finite element methods (FEMs) for curvature-driven flows has progressed more slowly. In particular, the convergence of various parametric FEMs for curvature flows of \emph{one-dimensional curves} and non-parametric FEMs for curvature flows of \emph{graph surfaces} {has been established in \cite{DeckelnickDziuk,Deckelnick-Dziuk-2009,Dziuk94,MR1681066,elliott2017approximations,Pozzi-Stinner-2021,Pozzi-Stinner-2023} and \cite{deckelnick1995convergence,Deckelnick-Dziuk-2006}, respectively.} The convergence of parametric FEMs for the mean curvature flow, Willmore flow, and surface diffusion of \emph{general closed surfaces} was established by Kov\'acs, Li, and Lubich~\cite{Kovacs-Li-Lubich-2019,KLL-Willmore} for finite elements of degree \(k \geq 2\) based on reformulating the governing equations using the evolution of the normal vector and mean curvature, with error estimates derived by comparing particle trajectories between the exact and approximate surfaces. More recently, the convergence of Dziuk's semi-implicit parametric FEM for the mean curvature flow of \emph{general closed surfaces} was proved in \cite{bai2023Dziuk} for finite elements of degree \(k \geq 3\), using a novel approach that estimates the distance between the exact and approximate surfaces. This distance-based error analysis, which neglects the tangential motion, also enabled the convergence proof of a stabilized BGN method for the curve shortening flow in~\cite{Bai-Li-MCOM2025}. Although the stabilized BGN method in \cite{Bai-Li-MCOM2025} bridges the BGN and MDR approaches and facilitates a convergence proof (at least for one-dimensional curves), it sacrifices the energy stability that is a hallmark of the original BGN method. The convergence of the original BGN method—remarkable for being both energy stable and capable of improving mesh quality through artificial tangential motion—remains an open and challenging problem.

Overall, the BGN method remains one of the most effective numerical approaches for simulating surface evolution under various curvature flows, due to its energy stability and mesh-improving properties. Its main drawback lies in a potential instability related to the choice of the time step size, which typically must be selected sufficiently large on a case-by-case basis in numerical simulations. This requirement may, however, increase the time discretization error.

The present paper aims to eliminate this drawback by bridging the BGN and MDR methods while preserving the energy stability of the original BGN approach. To this end, we propose a novel numerical scheme that not only ensures energy stability but also maintains high-quality meshes, even for small time step sizes. In order to clearly delineate the commonalities and differences between the BGN and MDR methods, we first consider their formulations when restricted to a properly defined tangential-motion space. More precisely, we define
\begin{align}\label{def-VTm-1}
\bar{V}^{(T),m-1} := \Big\{ \eta_h \in (S_h^{m-1})^d \,\Bigm|\, \int_{\Gamma_h^{m-1}}^{(h)} \eta_h \cdot n_h^{m-1}\,\chi_h = 0 \quad \forall\, \chi_h \in S_h^{m-1} \Big\}, 
\end{align}
which comprises those vector fields that are, in a mass-lumping quadrature sense, orthogonal to the normal direction \(n_h^{m-1}\). Within this tangential space, the BGN method in \eqref{BGN-weak-eqn} can be recast as
\[
\int_{\Gamma_h^{m-1}} \nabla_{\Gamma_h^{m-1}}\big(\tau v_h^m + \operatorname{id}\big) \cdot \nabla_{\Gamma_h^{m-1}} \eta_h = 0, \quad \forall\, \eta_h \in \bar{V}^{(T),m-1},
\]
while the finite element discretization of the MDR formulation \eqref{MDR-eqn} can be written as follows, by eliminating the term $ \int_{\Gamma_h^{m-1}}^{(h)} \kappa_h^m n_h^{m-1}\cdot \eta_h  $ using the definition of $\bar{V}^{(T),m-1}$ in \eqref{def-VTm-1}: 
\[
\int_{\Gamma_h^{m-1}} \nabla_{\Gamma_h^{m-1}} v_h^m \cdot \nabla_{\Gamma_h^{m-1}} \eta_h = 0, \quad \forall\, \eta_h \in \bar{V}^{(T),m-1}.
\]
The two formulations coincide when the test functions are further restricted to the subspace
\[
V^{(T),m-1} := \Big\{ \eta_h \in \bar{V}^{(T),m-1} \,\Bigm|\, \int_{\Gamma_h^{m-1}} \nabla_{\Gamma_h^{m-1}} \operatorname{id} \cdot \nabla_{\Gamma_h^{m-1}} \eta_h = 0 \Big\}.
\] 
Therefore, the velocities of the BGN and MDR methods agree on testing functions in $V^{(T),m-1}$ with the $H^1$ inner product. The difference between the velocities in this two methods lies in testing functions in the orthogonal complement of \(V^{(T),m-1}\) in \((S_h^{m-1})^d\), which is defined as \(V^{(N),m-1}\).

Let us denote by \(-\Delta_{\Gamma_h^{m-1}}: H^1(\Gamma_h^{m-1}) \to S_h^{m-1}\) the discrete Laplace-Beltrami operator, defined via a mass-lumping quadrature that fulfills
\[
\int_{\Gamma_h^{m-1}}^{(h)} \big(-\Delta_{\Gamma_h^{m-1}} f\big) \cdot \eta_h = \int_{\Gamma_h^{m-1}} \nabla_{\Gamma_h^{m-1}} f \cdot \nabla_{\Gamma_h^{m-1}} \eta_h,
\]
for any \(\eta_h \in S_h^{m-1}\) and for a fixed function \(f \in H^1(\Gamma_h^{m-1})\). For the vector-valued case, we define the operator
\(
-\Delta_{\Gamma_h^{m-1}}
\)
to act componentwise. {Then \(V^{(N),m-1}\) is the vector space given below:}
\begin{align*}
{V^{(N),m-1}
= \big\{\, \alpha \,(-\Delta_{\Gamma_h^{m-1}}{\rm id}) + I_h(n_h^{m-1}\chi_h)
:\ \alpha\in\mathbb{R},\ \chi_h\in S_h^{m-1} \,\big\}.}
\end{align*}
where \(I_h\) denotes the Lagrange interpolation operator {defined in \eqref{eq:def_Ih}}. Define \(T_h^{m-1} \in (S_h^{m-1})^d\) as the deviation from the space \(\mathrm{span}\left\{ I_h(n_h^{m-1}\,\chi_h) : \chi_h \in S_h^{m-1} \right\}\), i.e.,
\begin{align}\label{discrete-laplace}
    -\Delta_{\Gamma_h^{m-1}} {\rm{id}} = { I_h(\mu_h^{m-1}n_h^{m-1}) }  + T_h^{m-1},
\end{align}
with \(T_h^{m-1}\) orthogonal to the normal vector \(n_h^{m-1}\) at every finite element node. {From \eqref{discrete-laplace}, it follows that $\mu_h^{m-1}$ serves as an approximation of the curvature of $\Gamma_h^{m-1}$ at time level \(m-1\).} {Consequently, \(V^{(N),m-1}\) admits the following orthogonal decomposition:}
\begin{align*}
V^{(N),m-1}
&= \mathrm{span}\left\{ I_h(n_h^{m-1}\,\chi_h) : \chi_h \in S_h^{m-1} \right\} \oplus \mathrm{span}\left\{ T_h^{m-1} \right\} .
\end{align*}
Since the normal velocity of the surface is approximately determined by the specific curvature flow, the main difference between the velocities in the BGN and MDR methods lies in testing with $T_h^{m-1}$. 

This observation motivates the development of a numerical scheme that bridges the BGN and MDR methods by specifying the component of the velocity in the direction of \(T_h^{m-1}\), with the property of preserving energy stability while ensuring robust mesh quality with respect to the choice of the time step size. It is important to note that modifying the velocity \(v\) along the direction of \(T_h^{m-1}\) does not alter the fundamental velocity law governing the curvature flow, since \(T_h^{m-1}\) is a tangential vector that is orthogonal to the normal vector at every finite element node. Rather, this modification is introduced to achieve a more uniform distribution of mesh points, thereby substantially improving the quality of the computed surfaces.

The additional vector \(T_h^{m-1}\) defined in \eqref{discrete-laplace} quantifies the deviation of the polygonal curve or polyhedral surface—obtained from the fully discretized BGN scheme—from the conformal polyhedral surface or conformal polygonal curve introduced in \cite{BarrettGarckeNurnberg2020}. Conformal polygonal curves inherently yield an optimal mesh because, when the spatially semidiscrete or nonlinearly implicit BGN method for the evolution of one-dimensional curves admits a solution, adjacent edges are of equal length provided they are not parallel, see the discussion in \cite{BarrettGarckeNurnberg2020}. In particular, \(T_h^{m-1}\) emerges from the spatial discretization and the element-wise integration by parts performed on each curved triangle, thereby capturing the discontinuity in the conormal vector across adjacent curved triangles. In the case of planar curves discretized using the lowest-order FEM, if the two segments adjoining a node \(\zeta\) have unequal lengths, the vector \(T_h^{m-1}\) naturally points in the direction of the shorter segment. Accordingly, we prescribe a tangential velocity in the direction \(-T_h^{m-1}\) by enforcing the following constraint equation:
\begin{align}\label{add_constraint}
\int_{\Gamma_h^{m-1}}^{(h)} v_h^m \cdot T_h^{m-1} = - \alpha c^m\|T_h^{m-1}\|_{L^2_h}.
\end{align}
Here, \(\|\cdot\|_{L^2_h}\) denotes the discrete \(L^2\)-norm, and \(\alpha\) is a positive, adjustable parameter that can be varied across numerical experiments; unless specified otherwise, we set \(\alpha = 1\). The scalar variable \(\alpha c^m \in \mathbb{R}\), introduced via the constraint~\eqref{add_constraint}, is proportional to the magnitude of the velocity component in the \(-T_h^{m-1}\) direction and will be integrated into the formulation of the BGN-MDR scheme for both closed surfaces and closed curves.


In addition to simulating the evolution of closed surfaces under curvature flows, BGN methods have also been successfully applied to interface evolution problems, particularly in modeling solid-state dewetting processes involving contact line migration. Bao et al.\ have developed energy-stable parametric finite element methods that incorporate artificial tangential velocities in the spirit of the BGN framework~\cite{bao2022volume,bao2017parametric,Bao2021,Bao2023}. In this paper, we extend these ideas by defining the normal motion space and the auxiliary vector \(T_h^{m-1}\) for open surfaces with moving contact lines, as well as for open curves with moving contact points, in a manner analogous to the treatment of closed surfaces and curves. This extension properly incorporates the relevant boundary and contact angle conditions, and enables us to propose a BGN-MDR numerical scheme that ensures both energy stability and good mesh quality, even when small time step sizes are used.

The paper is organized as follows. In Section~2, we introduce the proposed BGN-MDR formulation for both mean curvature flow and surface diffusion on closed surfaces and closed curves in dimensions \(d=2,3\). A series of numerical experiments are presented to demonstrate that the method not only preserves energy stability but also maintains mesh quality even for very small time step sizes, in direct comparison with the BGN method. In Section~3, we extend the BGN-MDR scheme to the setting of open surfaces with a moving contact line in three dimensions, as well as to open curves with moving contact points in two dimensions. Several numerical examples are provided to verify the optimal temporal and spatial convergence rates and to illustrate the robustness of the proposed method in sustaining high-quality meshes for small time step sizes.

\section{BGN-MDR scheme on closed surfaces}\label{section:2}
Let \(\Gamma^0 = \Gamma(0)\subset\mathbb{R}^d\) denote the $(d-1)$-dimensional exact surface at the initial time \(t = 0\). We denote by \(\Gamma^0_h\) a piecewise triangular surface, where each element is the image of a reference triangle under a linear map, providing an approximation to \(\Gamma^0\).

Let \(t_m = m \tau\), \(m = 0, 1, \ldots, N\), be a partition of the time interval $[0,T]$, where \(\tau > 0\) is the time step size. For \(m \geq 1\), let \(x^{m-1}_j\), \(j = 1, \ldots, J\), denote the nodes of the numerical computed approximate surface \(\Gamma^{m-1}_h\) at time level \(t_{m-1}\).

Let \(\mathcal{K}^{m-1}_h\) denote the collection of triangles that constitute the approximate surface \(\Gamma^{m-1}_h\). The finite element space on \(\Gamma^{m-1}_h\) is defined by
\[
S_h^{m-1} = \left\{ v_h \in { C^0}(\Gamma^{m-1}_h) : v_h\big|_{K} \text{ is linear for all } K \in \mathcal{K}^{m-1}_h \right\}.
\]
In the vector-valued case, the finite element space \((S_h^{m-1})^l\), for any positive integer \(l\), is defined componentwise. The vector-valued finite element function \(X_h^m\in (S_h^{m-1})^d\) denotes the piecewise linear map from \(\Gamma_h^{m-1}\) onto \(\Gamma_h^m\), uniquely determined by its value at the finite element nodes.

\subsection{Tangential-motion space and normal-motion space}


We denote by $n_h^{m-1}$ the piecewise normal vector on the triangulated surface $\Gamma_h^{m-1}$, i.e., 
the restriction \(n_h^{m-1}|_K\) to each triangle \(K \subset \Gamma_h^{m-1}\) coincides with \(n_K^{m-1}\), the unit normal vector on \(K\). The discrete tangential-motion space $V^{(T),m-1}$ is defined as follows:
\begin{align}\label{eq:VT} 
    &V^{(T),m-1} := \\
    &\Big\{ \eta_h \in (S_h^{m-1})^d \,\big|\, (\eta_h, {n}_h^{m-1} \chi_h)^{(h)} = 0, \; \forall \chi_h \in S_h^{m-1} \text{ and } \int_{\Gamma_h^{m-1}} \nabla_{\Gamma_h^{m-1}} {\rm id} \cdot \nabla_{\Gamma_h^{m-1}} \eta_h = 0 \Big\}. \notag
\end{align} 
The space \(V^{(N),m-1}\) is defined as the orthogonal complement of \(V^{(T),m-1}\) in \((S_h^{m-1})^d\):
\begin{equation}\label{eq:VN}
V^{(N),m-1} = \big( V^{(T),m-1} \big)^\perp = {\rm{span}}\{ I_h(n_h^{m-1} \chi_h) \mid \chi_h \in S_h^{m-1}\} \oplus {\rm{span}}\{T_h^{m-1}\}.
\end{equation}
Here, the additional vector \(T_h^{m-1}\) is defined through the following problem: find
\[
(T_h^{m-1}, \mu_h^{m-1}) \in (S_h^{m-1})^d \times S_h^{m-1}
\]
such that
\begin{equation}\label{eq:T_h^m} 
\begin{aligned}
(T_h^{m-1}, \eta_h)^{(h)} + (\mu_h^{m-1}{n}_h^{m-1},  \eta_h)^{(h)} &= \int_{\Gamma_h^{m-1}} \nabla_{\Gamma_h^{m-1}} {\rm id} \cdot \nabla_{\Gamma_h^{m-1}} \eta_h
&& \forall \eta_h \in (S_h^{m-1})^d \\ 
(T_h^{m-1}, {n}_h^{m-1} \chi_h)^{(h)} &= 0 && \forall \chi_h \in S_h^{m-1} ,
\end{aligned} 
\end{equation} 
where \((\cdot,\cdot)^{(h)}\) denotes the mass-lumping inner product on the discrete surface \(\Gamma_h^{m-1}\). {Denote by $\{\sigma_l\}_{l=1}^L$ a family of pairwise disjoint, relatively open $(d-1)$-simplices that constitute $\Gamma_h^{m-1}$. For piecewise continuous functions \(\eta, \chi \in L^{\infty}(\Gamma_h^{m-1})\), which may have jumps across the edges of \(\{\sigma_l\}_{l=1}^L\), the mass-lumping inner product is defined as follows (see Definition 43 in \cite{BarrettGarckeNurnberg2020}):
\begin{align}\label{mass-lump-def}
( \eta, \chi )^{(h)} = \frac{1}{d} \sum_{l=1}^L \mathcal{H}^{d-1}(\sigma_l) \sum_{k=1}^{d} \left(\eta \cdot \chi\right)\left((\vec{q}_{l,k})^{-}\right),
\end{align}
where $\vec{q}_{l,k}$ denotes the $k$-th vertex of the $l$-th simplex, and $\mathcal{H}^{d-1}$ is the $d-1$-dimensional Hausdorff measure; here, $\eta\left((\vec{q}_{l,k})^{-}\right) = \lim_{\sigma_l \ni \vec{p} \to \vec{q}_{l,k}} \eta(\vec{p})$. This definition is naturally extended to vector- and tensor-valued functions, with the product between $\eta$ and $\chi$ understood as the (Euclidean) inner product.}
{For any function $f$ defined on $\Gamma_h^{m-1}$ , possibly piecewise defined and discontinuous across the boundaries of the triangles, we define its Lagrange interpolation $I_hf \in S_h^{m-1}$ as follows:
\begin{equation}\label{eq:def_Ih}
    (I_hf, w_h)^{(h)} = (f, w_h)^{(h)}, \quad \forall w_h \in S_h^{m-1}.
\end{equation}}


By introducing an averaged normal vector at each finite element node, which incorporates the contributions from all neighboring simplices weighted by their respective areas (in three dimensions) or lengths (in two dimensions), we define
\begin{align}\label{averaged-normal}
    \hat{n}_h^{m-1}(x_j^{m-1}) := \frac{\sum_{\sigma_l \ni x_j^{m-1}} |\sigma_l|\, n_{\sigma_l}^{m-1}}{\left|\sum_{\sigma_l \ni x_j^{m-1}} |\sigma_l|\, n_{\sigma_l}^{m-1}\right|},
\end{align}
{where \(|\sigma_l| := \mathcal{H}^{d-1}(\sigma_l)\) denotes the area (for \(d=3\)) or length (for \(d=2\)) of the simplex \(\sigma_l\), and \(n_{\sigma_l}^{m-1}\) is the constant unit normal vector associated with \(\sigma_l\). We use the same notation, \(\hat n_h^{m-1}\), to denote the finite element function whose value at each node \(x_j^{m-1}\) coincides with the averaged normal vector \(\hat n_h^{m-1}(x_j^{m-1})\).}
{Furthermore, from the mass lumping definition in \eqref{mass-lump-def}, we obtain
\begin{align*}
(\mu_h^{m-1}{n}_h^{m-1},  \eta_h)^{(h)}
&= \frac{1}{d} \sum_{l=1}^L |\sigma_l| \sum_{k=1}^{d} \left(\mu_{h}^{m-1} n_{\sigma_l}^{m-1} \cdot \eta_h\right)\left((\vec{q}_{l,k})^{-}\right) \\
&= \frac{1}{d} \sum_{j=1}^J \big(\mu_h^{m-1}(x_j^{m-1})\,\eta_h(x_j^{m-1}) \cdot \sum_{\sigma_l \ni x_j^{m-1}} |\sigma_l|\, n_{\sigma_l}^{m-1}\big) \\
&= \frac{1}{d} \sum_{j=1}^J \Big(\Big| \sum_{\sigma_l \ni x_j^{m-1}} |\sigma_l|\, n_{\sigma_l}^{m-1} \Big|
\mu_h^{m-1}(x^{m-1}_j)\, \hat n_h^{m-1}(x^{m-1}_j) \cdot \eta_h(x^{m-1}_j)\Big).
\end{align*}}

The system \eqref{eq:T_h^m} can be shown to be well-posed by demonstrating that its associated homogeneous linear system admits only the trivial solution. Specifically, choosing the test functions \(\eta_h := T_h^{m-1}\) and \(\chi_h := \mu_h^{m-1}\) in \eqref{eq:T_h^m} leads to \(T_h^{m-1} \equiv 0\). Subsequently, selecting the test function {\(\eta_h := I_h(\hat{n}_h^{m-1} \mu_h^{m-1} )\)} yields \(\mu_h^{m-1} \equiv 0\), thereby establishing the desired well-posedness. 



\subsection{Numerical scheme for mean curvature flow}
The proposed BGN-MDR scheme for mean curvature flow is formulated as follows: Determine 
\[
(v_h^m, \lambda_h^m,  c^m) \in (S_h^{m-1})^d \times S_h^{m-1} \times \mathbb{R}
\]
such that
\begin{subequations}\label{eq:MCF-num-equ}
    \begin{align}
        \int_{\Gamma_h^{m-1}} \nabla_{\Gamma_h^{m-1}} (\tau v_h^m + {\rm id}) \cdot \nabla_{\Gamma_h^{m-1}} \eta_h - \left( ({n}_h^{m-1}\lambda_h^m, \eta_h)^{(h)} + c^m (T_h^{m-1}  , \eta_h)^{(h)}\right) &= 0 \label{eq:MCF-num-1-equ} \\ 
        (v_h^m \cdot {n}_h^{m-1}, \chi_h)^{(h)} + (\lambda_h^m, \chi_h)^{(h)} &= 0  \label{eq:MCF-num-2-equ} \\[5pt]
        (v_h^m , T_h^{m-1})^{(h)} + \alpha c^m \|T_h^{m-1}\|_{L^2_h}  &= 0 \label{eq:MCF-num-3-equ}
    \end{align}    
\end{subequations}
for all $\eta_h \in (S_h^{m-1})^d$ and $\chi_h \in S_h^{m-1}$. {Here, $\lambda_h^m$ serves as an approximation of the curvature of $\Gamma_h^m$ at time level $m$, which is distinct from $\mu_h^{m-1}$ in \eqref{discrete-laplace}, as the latter approximates the curvature of the previous time level surface
}. In practical computations, the normalized term \(T_h^{m-1} / \|T_h^{m-1}\|_{L^2_h}\) is homogeneous of degree zero with respect to the vector \(T_h^{m-1}\), and as a result, it remains uniformly bounded as \(\|T_h^{m-1}\|_{L^2_h}\) approaches the order of machine precision. In such cases, the triangulated surface forms a conformal polyhedron with good mesh quality (see~\cite{BarrettGarckeNurnberg2020}), and the definition of \((v_h^m, T_h^{m-1} / \|T_h^{m-1}\|_{L^2_h})^{(h)}\) will only drive the mesh slightly away from being conformal polyhedral and therefore does not immediately cause mesh distortion. Based on this observation, it is not necessary to treat the degenerate case \(\|T_h^{m-1}\|_{L^2_h} = 0\) separately. Therefore, for the remainder of our analysis on the effectiveness of the numerical scheme \eqref{eq:MCF-num-equ}, we shall assume \(\|T_h^{m-1}\|_{L^2_h} \neq 0\), without further clarification. 

From equations \eqref{eq:MCF-num-2-equ} and \eqref{eq:MCF-num-3-equ}, \( \lambda_h^m \in S_h^{m-1} \) and \(c^m \in \mathbb{R} \) can be solved as follows:
\[
\lambda_h^m = -I_h(v_h^m \cdot {n}_h^{m-1}) \quad \text{and} \quad 
 c^m = - \Big(v_h^m , \frac{T_h^{m-1}}{\alpha\|T_h^{m-1}\|_{L^2_h}} \Big)^{(h)}.
\]
Substituting the above expressions into \eqref{eq:MCF-num-1-equ} yields:
\begin{align}\label{eq:MCF-scheme--equal}
          \quad\int_{\Gamma_h^{m-1}} \nabla_{\Gamma_h^{m-1}} (\tau v_h^m + {\rm id}) \cdot \,\nabla_{\Gamma_h^{m-1}} \eta_h &+ (v_h^m \cdot {n}_h^{m-1} , {n}_h^{m-1} \cdot \eta_h)^{(h)} \\
    &\hspace{-20pt}+\frac{(v_h^m , T_h^{m-1})^{(h)} (T_h^{m-1}, \eta_h)^{(h)}}{\alpha\|T_h^{m-1}\|_{L^2_h}} = 0
    \quad\forall\, \eta_h \in (S_h^{m-1})^d. \notag
\end{align} 
This formulation shows that the BGN-MDR method reduces to the BGN method as $\alpha\rightarrow\infty$. 

Under some mild conditions stated in Theorem \ref{thm:well-pose-mass-mcf}, the well-posedness of the numerical scheme \eqref{eq:MCF-num-equ} will be established.

\begin{theorem}\label{thm:well-pose-mass-mcf}
    Assume that the discrete surface (or the discrete curve) $\Gamma_h^{m-1}$ satisfies the following mild conditions:
    \begin{enumerate}[label={(A\arabic*)}]
        \item The elements are nondegenerate, i.e., for each $K\subset \mathcal{K}_h^{m-1}$, it holds
        $ |K| > 0$.

        \item The normal vectors 
        fulfill
        \[
        \dim \bigl(\operatorname{span}\big\{(n_h^{ m-1}\phi_h,1)^{(h)}\big|\,\phi_h \in S_h^{m-1}\big\} \bigr) \;=\; d.
        \]
        \end{enumerate}
        
        \noindent
    Then, the numerical scheme \eqref{eq:MCF-num-equ} is well-posed, i.e., there exists a unique solution \( (v_h^m, \lambda_h^m,  c^m) \in (S_h^{m-1})^d \times S_h^{m-1} \times \mathbb{R} \).
\end{theorem}

\begin{proof}
It suffices to demonstrate that the following homogeneous system admits only the trivial solution: 
\begin{subequations}\label{well-pose-homo}
    \begin{align}
        \int_{\Gamma_h^{m-1}} \nabla_{\Gamma_h^{m-1}} \tau v_h^{m} \cdot \nabla_{\Gamma_h^{m-1}} \eta_h - \Bigl( \bigl(n_h^{m-1}\lambda_h^m, \eta_h\bigr)^{(h)} +  c^m \bigl(T_h^{m-1}  , \eta_h\bigr)^{(h)}\Bigr) &= 0 
        \label{well-pose-homo1} \\[1mm]
        \bigl(v_h^m \cdot n_h^{m-1}, \chi_h\bigr)^{(h)} + \bigl(\lambda_h^m, \chi_h\bigr)^{(h)} &= 0 \label{well-pose-homo2} \\[1mm]
        \bigl(v_h^m , T_h^{m-1}\bigr)^{(h)} + \alpha c^m \|T_h^{m-1}\|_{L^2_h}  &= 0  \label{well-pose-homo3}
    \end{align}    
\end{subequations}
for $\eta_h \in (S_h^{m-1})^d$ and $\chi_h \in S_h^{m-1}$. 
By choosing the test functions \(\eta_h := v_h^m\) and \(\chi_h := \lambda_h^m\) in \eqref{well-pose-homo}, solving for \( c^m\) from \eqref{well-pose-homo3}, and then summing the equations \eqref{well-pose-homo1} and \eqref{well-pose-homo2}, we obtain
\begin{align}\label{well-pose-homo-1}
    \tau \|\nabla_{\Gamma_h^{m-1}} v_h^m\|_{L^2(\Gamma_h^{m-1})}^2 + \|v_h^m \cdot n_h^{m-1}\|_{L^2_h}^2 
    + \|T_h^{m-1}\|_{L^2_h}^{(h)}\Big|\Big(v_h^m , \frac{T_h^{m-1}}{\alpha\,\|T_h^{m-1}\|_{L^2_h}}\Big)^{(h)}\Big|^2 = 0,
\end{align}
which implies 
\begin{subequations}\label{well-pose-homo-2}
    \begin{align}
        \|\nabla_{\Gamma_h^{m-1}} v_h^m\|_{L^2(\Gamma_h^{m-1})} &= 0, \label{well-pose-homo-2-a}\\[1mm]
        I_h(v_h^m \cdot n_h^{m-1}) &= 0. \label{well-pose-homo-2-b}
    \end{align}
\end{subequations}
From the nondegenerate condition \((A1)\), equation \eqref{well-pose-homo-2-a} immediately implies that \(v_h^m\) is a constant vector field, that is, there exists a constant vector \(v_{\rm constant}\) such that \(v_h^m = v_{\rm constant}\). In conjunction with \eqref{well-pose-homo-2-b}, this entails
\[
v_{\rm constant} \cdot (n_h^{ m-1} \phi_h, 1)^{(h)} = 0.
\]
Hence, by 
condition \((A2)\), one deduces that \(v_{\rm constant} = 0\), implying \(v_h^m \equiv 0\). Substituting \(v_h^m = 0\) into \eqref{well-pose-homo2} immediately yields \(\lambda_h^m = 0\), and then from \eqref{well-pose-homo3} it follows that
\[
c^m = -\Big(v_h^m, \frac{T_h^{m-1}}{\alpha\,\|T_h^{m-1}\|_{L^2_h}}\Big)^{(h)} = 0.
\]
\hfill 
\end{proof}

{ 
\begin{remark}\upshape 
Condition \((A2)\) requires that the discrete vertex normals of \(\Gamma_h^{m-1}\) span \(\mathbb{R}^d\). This condition is violated only in exceptional cases; for example, it is always satisfied for surfaces \(\Gamma_h^{m-1}\) without self-intersections (see Remark 65 in \cite{BarrettGarckeNurnberg2020}).
This mild assumption~$(A2)$ is widely used in proving the well-posedness of numerical schemes for curvature flows; see, for example,~\cite{Bao2023,li2021energy,garcke2025stable}.
\end{remark} }

Moreover, the following theorem states that the proposed BGN-MDR scheme is unconditionally energy stable (area decreasing) for mean curvature flow.

\begin{theorem}\label{thm:area-decreasing-mass-mcf}  
Let \( v_h^m \in (S_h^{m-1})^d \) denote the solution of \eqref{eq:MCF-scheme--equal}. Then the area of the discrete surface satisfies
\[
|\Gamma_h^{m}| \le |\Gamma_h^{m-1}| \quad \text{for } m = 1, 2, \ldots, N,
\]
where \( |\Gamma_h^{m}| \) represents the area of the numerical surface \( \Gamma_h^{m} \).
\end{theorem}

\begin{proof}
Choosing the test function \( \eta_h = v_h^m \) in \eqref{eq:MCF-scheme--equal} leads to
\begin{align}\label{thm:area-decreasing-mass-mcf-2}
    \int_{\Gamma_h^{m-1}} \nabla_{\Gamma_h^{m-1}} (\tau v_h^m + {\rm id}) \cdot \nabla_{\Gamma_h^{m-1}} v_h^m 
    + \left(v_h^m \cdot n_h^{m-1},\, v_h^m \cdot n_h^{m-1}\right)^{(h)}
    + \frac{\left[(v_h^m, T_h^{m-1})^{(h)}\right]^2}{\alpha\|T_h^{m-1}\|_{L^2_h}}= 0.
\end{align}
By substituting the following inequality (see \cite[(2.21)]{barrett2008hypersurfaces} or \cite[(2.31)]{barrett2007parametric}) into \eqref{thm:area-decreasing-mass-mcf-2}
\begin{align}\label{area-decreasing-important-property}
\int_{\Gamma_h^{m-1}} \nabla_{\Gamma_h^{m-1}} (\tau v_h^m + {\rm id}) \cdot \nabla_{\Gamma_h^{m-1}} v_h^m  \ge |\Gamma_h^m| - |\Gamma_h^{m-1}| , 
\end{align}
we obtain 
$|\Gamma_h^m| - |\Gamma_h^{m-1}| \le 0$.
\hfill\end{proof}

We demonstrate the advantages of the BGN-MDR method in improving mesh quality through the following examples.

\begin{example}\label{Example1}\upshape
We consider a torus-shape surface in Figure \ref{Example1a}, defined by the parametrization
\[
    x = \begin{bmatrix}
    (1 + 0.65 \cos \varphi) \cos \theta \\ 
    (1 + 0.65 \cos \varphi) \sin \theta \\
    0.65 \sin \varphi + 0.3 \sin(5 \theta)
    \end{bmatrix},
    \quad \theta \in [0, 2\pi], \quad \varphi \in [0, 2\pi].
\]
Numerical simulations of mean curvature flow starting from the torus-shaped initial surface described above, obtained by using the BGN and BGN-MDR schemes with 5592 triangles, are shown in Figure~\ref{eg:Example1}. Both methods achieve good mesh quality with a time step size of \(\tau = 0.005\), as illustrated in Figures~\ref{Example1d} and \ref{Example1e}. However, when the step size is reduced to \(\tau = 10^{-4}\), the BGN scheme tends to exhibit mesh distortion, whereas the BGN-MDR scheme continues to preserve good mesh quality, as shown in Figures~\ref{Example1b} and \ref{Example1c}. 

Furthermore, Figure~\ref{eg:Example1-auxiliary} shows the energy stability (area decrease) of the BGN-MDR scheme and the evolution of $\|T_h^{m-1}\|_{L^2_h}$. These numerical results demonstrate the superiority of the BGN-MDR method over the BGN method in improving mesh quality, while maintaining energy stability. 
    

    \begin{figure}[htbp]
        \centering
        
        \begin{subfigure}[b]{0.25\textwidth}
        \centering
        \includegraphics[width=0.9\textwidth]{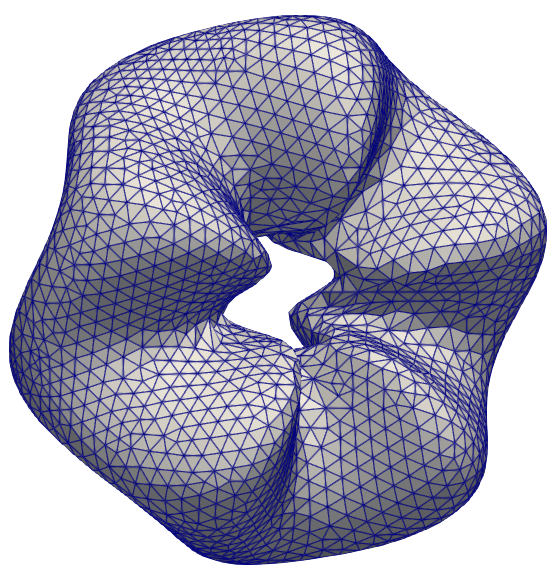}
        \caption{Initial Surface}
        \label{Example1a}
        \end{subfigure}
        \vspace{5pt} 
        
        \begin{subfigure}[b]{0.45\textwidth}
        \centering
            \includegraphics[width=0.45\textwidth]{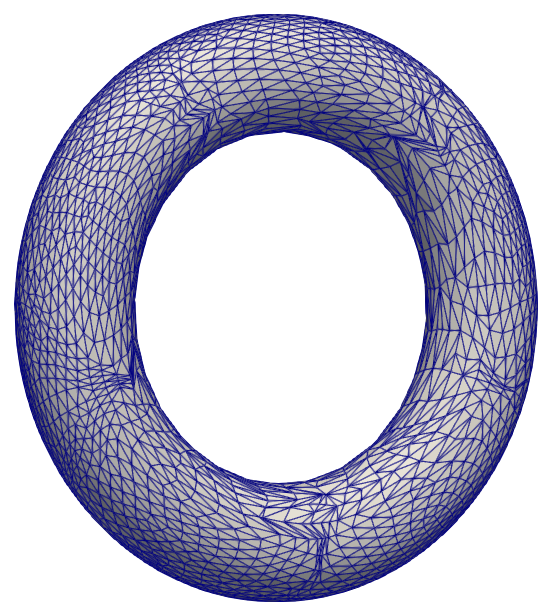}
            \caption{BGN scheme for $\tau=10^{-4}$}
            \label{Example1b}
        \end{subfigure}
        \hspace{10pt}
        \begin{subfigure}[b]{0.45\textwidth}
        \centering
            \includegraphics[width=0.45\textwidth]{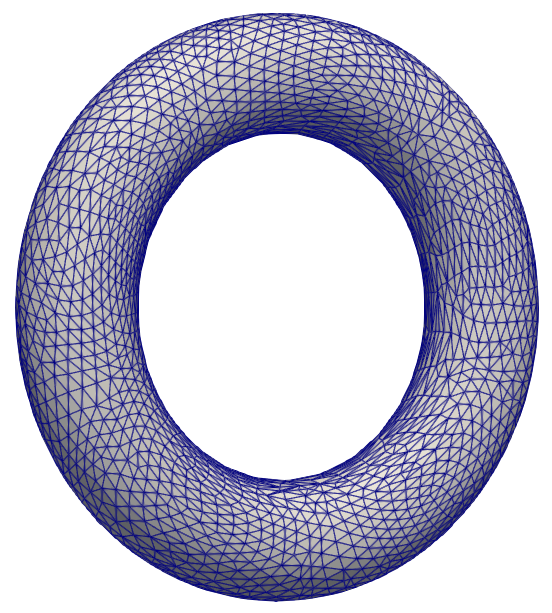}
            \caption{BGN-MDR scheme for $\tau=10^{-4}$}
            \label{Example1c}
        \end{subfigure}
        
        \vspace{5pt} 
        
        \begin{subfigure}[b]{0.45\textwidth}
        \centering
            \includegraphics[width=0.45\textwidth]{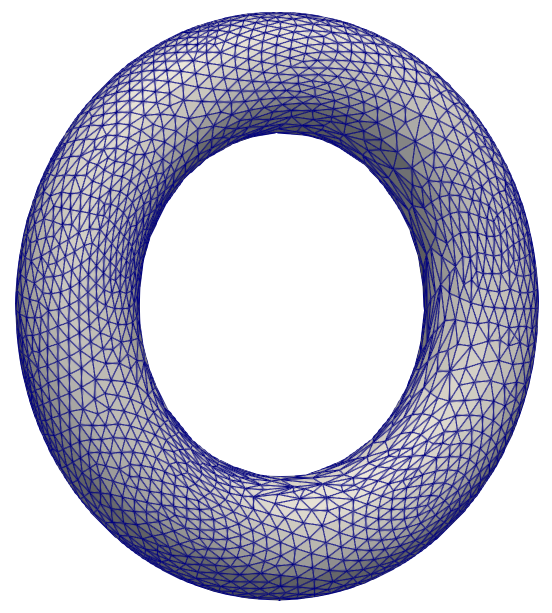}
            \caption{BGN scheme for $\tau=0.005$}
            \label{Example1d}
        \end{subfigure}
        \hspace{10pt}
        \begin{subfigure}[b]{0.45\textwidth}
        \centering
            \includegraphics[width=0.45\textwidth]{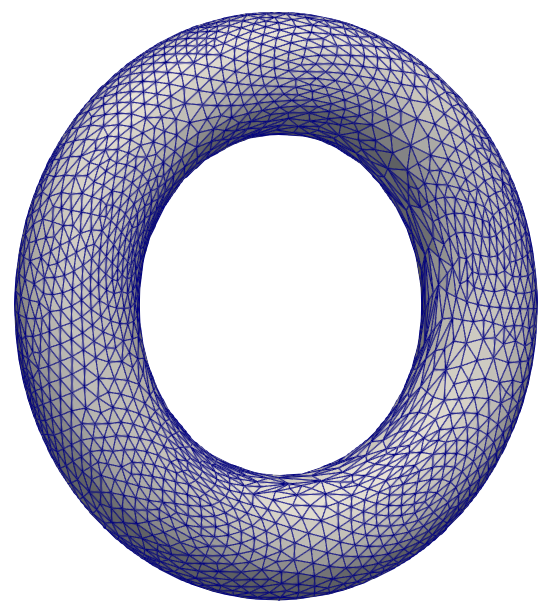}
            \caption{BGN-MDR scheme for $\tau=0.005$}
            \label{Example1e}
        \end{subfigure}
        
        \vspace{-5pt}
        \caption{Surface evolution in Example \ref{Example1}}
        \label{eg:Example1}
    \end{figure}
    
    \begin{figure}[htbp]
        \centering
        \begin{subfigure}[b]{0.45\textwidth}
            \includegraphics[width=0.9\textwidth]{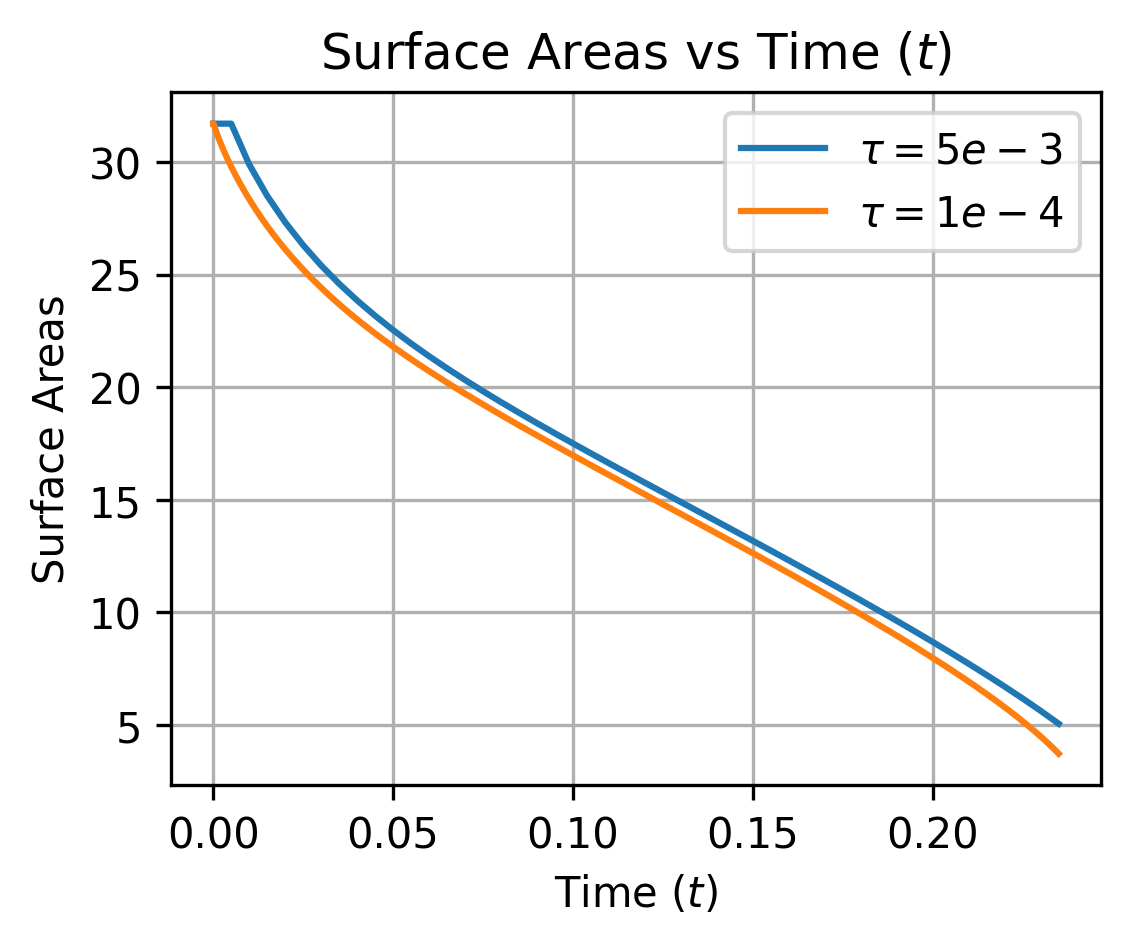}
            \caption{Evolution of surface area}
            \label{Example1-area-1e-4}
        \end{subfigure}
    \hfill
%
        \begin{subfigure}[b]{0.45\textwidth}
            \includegraphics[width=0.9\textwidth]{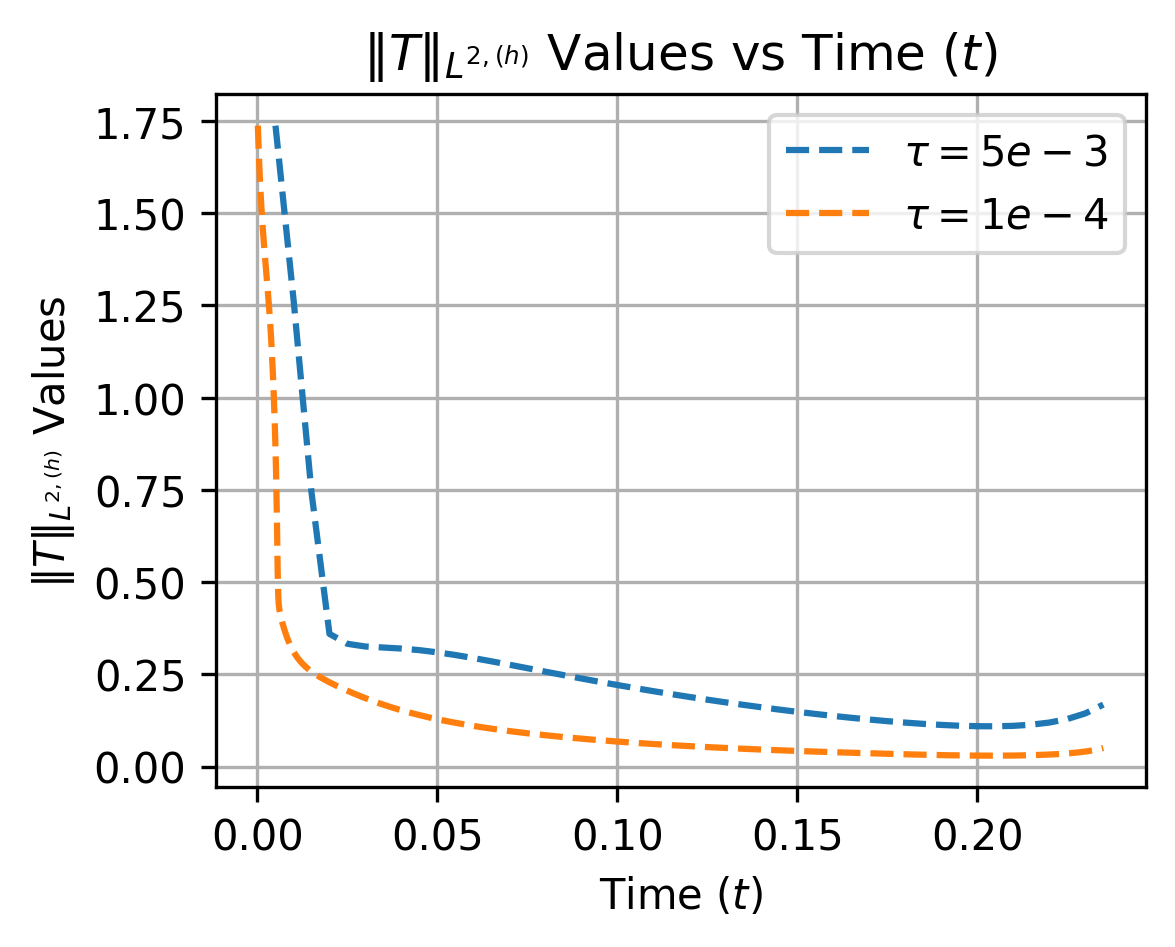}
            \caption{Evolution of $\|T\|_{L^2_h}$}
            \label{Example1-T2-1e-4}
        \end{subfigure}
    
        \vspace{-5pt}
        \caption{Evolution of surface area and $\|T\|_{L^2_h}$ in Example \ref{Example1}.}
        \label{eg:Example1-auxiliary}
    \end{figure}

    \end{example}

    \begin{example}\label{Example3}\upshape
        We consider the dumbbell-shape surface shown in Figure \ref{Example3a}, defined by the parametrization
        \[
        x = \begin{bmatrix}
        (0.6 \cos^2 \varphi + 0.4) \cos \theta \sin \varphi \\ 
        (0.6 \cos^2 \varphi + 0.4) \sin \theta \sin \varphi \\
        \cos \varphi
        \end{bmatrix},
        \quad \theta \in [0, 2\pi], \quad \varphi \in [0, 2\pi].
        \]
Under mean curvature flow, the dumbbell-shaped surface evolves into a sphere which then shrinks to a point singularity. To enhance the stability of the BGN-MDR algorithm in the presence of the blow-up of \(\|T_h^{m-1}\|_{L^\infty}\) near the singularity, for \(t\geq 0.0908\) we modify the constraint equation \eqref{eq:MCF-num-3-equ} by incorporating the factor
\[
\frac{(\mathrm{area}(\Gamma_h^{m-1}))^{1/2}}{(\mathrm{area}|_{t= 0.0908})^{1/2}}  
\]
which compensates for the surface area diminishing to zero as the evolution progresses. The modified constraint equation can be written as follows:
\begin{align}\label{add_constraint_compen}
\left(v_h^m, \frac{-T_h^{m-1}}{\|T_h^{m-1}\|_{L^2_h}}\right)^{(h)} +  c^m \cdot \frac{(\mathrm{area}(\Gamma_h^{m-1}))^{1/2}}{(\mathrm{area}|_{t= 0.0908})^{1/2}} = 0.
\end{align}

We employ the BGN method and the proposed BGN-MDR method with mesh size \(h=0.06\), with 3836 triangles and 1920 vertices. To accurately resolve the solution near the blow-up time, the time step size is reduced from \(\tau_1=10^{-4}\) to \(\tau_2=2\times10^{-7}\) for \(t\geq 0.0908\). Both the BGN and BGN-MDR schemes yield good mesh quality under mean curvature flow; see Figures~\ref{Example3d} and \ref{Example3e}.

However, when the time-step is reduced from \(\tau_1=10^{-5}\) to \(\tau_2=10^{-7}\), the BGN scheme tends to exhibit deteriorated mesh quality, whereas the BGN-MDR scheme maintains a high-quality mesh; see Figures~\ref{Example3b} and \ref{Example3c}. 
%
These results demonstrate the superiority of the BGN-MDR method in improving the mesh quality.


        \begin{figure}[htbp]
            \centering
            
            \begin{subfigure}[b]{0.35\textwidth}
                \includegraphics[width=0.9\textwidth]{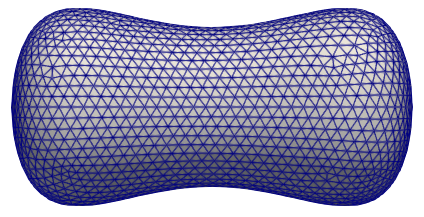}
                \caption{Initial Surface}
                \label{Example3a}
            \end{subfigure}
            
            \vspace{5pt} 
            
            \begin{subfigure}[b]{0.4\textwidth}
            \centering
                \includegraphics[width=0.45\textwidth]{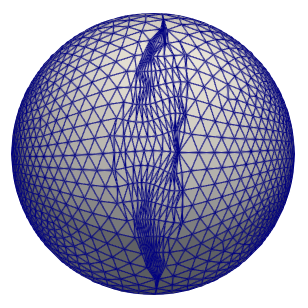}
                \caption{BGN with $\tau_1=10^{-5}$ and $\tau_2=10^{-7}$ at $T=0.0906655$}
                \label{Example3b}
            \end{subfigure}
            \hspace{20pt}
            \begin{subfigure}[b]{0.4\textwidth}
            \centering
                \includegraphics[width=0.45\textwidth]{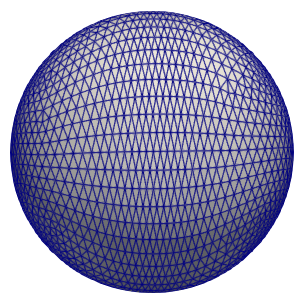}
                \caption{BGN-MDR with $\tau_1=10^{-5}$, $\tau_2=10^{-7}$ at $T=0.0906106$}
                \label{Example3c}
            \end{subfigure}
            
            \vspace{5pt} 
            
            \begin{subfigure}[b]{0.4\textwidth}
            \centering
                \includegraphics[width=0.45\textwidth]{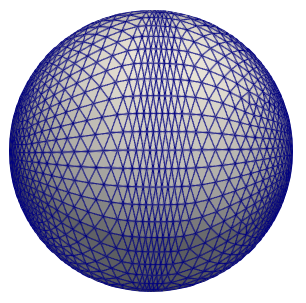}
                \caption{BGN with $\tau_1=10^{-4}$, $\tau_2=2\times 10^{-7}$ at $T=0.0908936$}
                \label{Example3d}
            \end{subfigure}
            \hspace{20pt}
            \begin{subfigure}[b]{0.4\textwidth}
            \centering
                \includegraphics[width=0.45\textwidth]{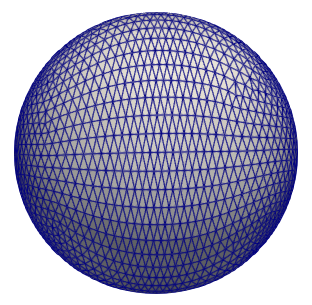}
                \caption{BGN-MDR with $\tau_1=10^{-4}$, $\tau_2=2\times 10^{-7}$ at $T=0.0909028$}
                \label{Example3e}
            \end{subfigure}
            
            \vspace{-5pt}
            \caption{Surface evolution in Example \ref{Example3}}
            \label{eg:Example3}
        \end{figure}
        
%

        \end{example}

\subsection{Numerical scheme for surface diffusion}
In the BGN-MDR scheme for mean curvature flow, the constraint equation \eqref{eq:MCF-num-3-equ} imposes an artificial tangential velocity along the direction \(T_h^{m-1}\), as defined in \eqref{eq:T_h^m}. By adopting this constraint, we formulate the BGN-MDR scheme for surface diffusion as follows: Find 
\[
(v_h^m, \lambda_h^m, c^m) \in (S_h^{m-1})^d \times S_h^{m-1} \times \mathbb{R}
\]
such that
\begin{subequations}\label{eq:SD-num}
\begin{align}
\int_{\Gamma_h^{m-1}} \nabla_{\Gamma_h^{m-1}} (\tau v_h^m + {\rm id}) \cdot \nabla_{\Gamma_h^{m-1}} \eta_h  - (\lambda_h^m\, {n}_h^{m-1}, \eta_h)^{(h)} - c^m (T_h^{m-1}, \eta_h)^{(h)} &= 0 \label{eq:SD-num-1} \\
(v_h^m \cdot n_h^{m-1}, \chi_h)^{(h)} + \int_{\Gamma_h^{m-1}} \nabla_{\Gamma_h^{m-1}} \lambda_h^m \cdot \nabla_{\Gamma_h^{m-1}} \chi_h &= 0 \label{eq:SD-num-2} \\
(v_h^m, T_h^{m-1})^{(h)} + \alpha c^m \|T_h^{m-1}\|_{L^2_h}  &= 0  \label{eq:SD-num-3}
\end{align}
\end{subequations}
for all $\eta_h \in (S_h^{m-1})^d$ and $\chi_h \in S_h^{m-1}$. 



The well-posedness of the numerical scheme in \eqref{eq:SD-num} is shown in the following theorem.
 
\begin{theorem}\label{thm:SD-num}
    Assume that the discrete surface (or the discrete curve) $\Gamma_h^{m-1}$ satisfies the following mild conditions:
    \begin{enumerate}[label={(A\arabic*)}]
        \item The elements are nondegenerate, i.e., for each $K\subset \mathcal{K}_h^{m-1}$, it holds
        $
        |K| > 0.
        $ 
        \item The normal vectors 
        fulfill 
        \[
        \dim \bigl(\operatorname{span}\big\{(n_h^{m-1}\phi_h,1)^{(h)}\big|\,\phi_h \in S_h^{m-1}\big\} \bigr) \;=\; d.
        \]
        \end{enumerate}
        
        \noindent
   Then the numerical scheme in \eqref{eq:SD-num} possesses a unique solution
    \[
    (v_h^m, \lambda_h^m, c^m) \in (S_h^{m-1})^d \times S_h^{m-1} \times \mathbb{R}.
    \]
\end{theorem}

\begin{proof}
    In order to establish the existence and uniqueness of solutions to \eqref{eq:SD-num}, it is sufficient to show that the following homogeneous system admits only the trivial solution:
    \begin{subequations}\label{homo-mass-SD}
        \begin{align}
            \int_{\Gamma_h^{m-1}} \nabla_{\Gamma_h^{m-1}} \tau v_h^{m} \cdot \nabla_{\Gamma_h^{m-1}} \eta_h - ( \lambda_h^m n_h^{m-1}, \eta_h)^{(h)} -  c^m (T_h^{m-1} , \eta_h)^{(h)} &= 0\label{homo-SD-1} \\ 
            (v_h^m \cdot {n}_h^{m-1}, \chi_h)^{(h)} + \int_{\Gamma_h^{m-1}} \nabla_{\Gamma_h^{m-1}} \lambda_h^m \cdot \nabla_{\Gamma_h^{m-1}} \chi_h &= 0 \label{homo-SD-2} \\ 
            (v_h^m , T_h^{m-1})^{(h)} + \alpha c^m \|T_h^{m-1}\|_{L^2_h} &= 0 \label{homo-SD-3}
        \end{align}    
    \end{subequations}
for all $\eta_h \in (S_h^{m-1})^d$ and $\chi_h \in S_h^{m-1}$. 

By choosing \( \eta_h := v_h^m \) in \eqref{homo-SD-1} and \( \chi_h := \lambda_h^m \) in \eqref{homo-SD-2}, and then summing equations \eqref{homo-SD-1}--\eqref{homo-SD-3}, we obtain
    \begin{align}\label{homo-mass-SD-temp}
        \tau\|\nabla_{\Gamma_h^{m-1}} v_h^{m}\|_{L^2(\Gamma_h^{m-1})}^2 + \|\nabla_{\Gamma_h^{m-1}} \lambda_h^m\|_{L^2(\Gamma_h^{m-1})}^2 + \frac{1}{\alpha}| c^m|^2/\|T_h^{m-1}\|_{L^2_h} = 0.
    \end{align}
    Consequently, it follows that
    \begin{subequations}\label{well-pose-mass-SD}
        \begin{align}
            \|\nabla_{\Gamma_h^{m-1}} v_h^{m}\|_{L^2(\Gamma_h^{m-1})} &= 0, \label{well-pose-mass-SD-1} \\
            \|\nabla_{\Gamma_h^{m-1}} \lambda_h^m\|_{L^2(\Gamma_h^{m-1})} &= 0, \label{well-pose-mass-SD-2} \\
            \frac{1}{\alpha}| c^m|^2\|T_h^{m-1}\|_{L^2_h} &= 0. \label{well-pose-mass-SD-3}
        \end{align}
    \end{subequations}
    Substituting \eqref{well-pose-mass-SD-2} into the weak formulation \eqref{homo-SD-2} yields
    \begin{align}\label{homo-temp-v}
        I_h (v_h^m \cdot n_h^{m-1}) = 0.
    \end{align}
From the nondegenerate condition \((A1)\) and the  condition \((A2)\), equations \eqref{well-pose-mass-SD-1} and \eqref{homo-temp-v} imply that \(v_h^m \equiv 0\). 
Furthermore, by selecting \( \eta_h := I_h (\lambda_h^m {n}_h^{m-1}) \) in \eqref{homo-SD-1} and noting that \( v_h^m \equiv 0 \) along with the condition \( (T_h^{m-1}, \lambda_h^m {n}_h^{m-1})^{(h)} = 0 \) inherent in the definition of \( T_h^{m-1} \), one obtains
    \begin{align}\label{homo-temp-lambda}
        (\lambda_h^m {n}_h^{m-1}, \lambda_h^m {n}_h^{m-1})^{(h)} = 0,
    \end{align}
    which implies that \( \lambda_h^m \equiv 0 \). Finally, employing \( v_h^m \equiv 0 \) in conjunction with the constraint \eqref{well-pose-mass-SD-3} immediately yields
    \begin{align}\label{homo-temp-c}
        c^m = \frac{1}{\alpha} \Big(v_h^m, \frac{-T_h^{m-1}}{\|T_h^{m-1}\|_{L^2_h}}\Big)^{(h)}= 0.
    \end{align}
Therefore, the homogeneous linear system in \eqref{homo-mass-SD} admits only the trivial solution, which in turn implies that the original linear system in \eqref{eq:SD-num} possesses a unique solution. 
\hfill\end{proof}



The following theorem states that the BGN-MDR scheme for surface diffusion is unconditionally energy stable (area decreasing).  

\begin{theorem}\label{thm:area-decreasing-mass-sd}  
    Let \((v_h^m, \lambda_h^m,  c^m) \in (S_h^{m-1})^d \times S_h^{m-1} \times \mathbb{R}\) be the numerical solution of surface diffusion determined by \eqref{eq:SD-num}. Then the following inequality holds: 
    \begin{align}\label{thm:area-decreasing-mass-sd-1}
        |\Gamma_h^{m}| \le |\Gamma_h^{m-1}| \quad \text{for } m = 1, 2, \ldots, N,
    \end{align}
where \( |\Gamma_h^{m}| \) denotes the area of the numerical surface \( \Gamma_h^{m} \).
\end{theorem}

\begin{proof}
By choosing the test function \( \eta_h := v_h^m  \) in \eqref{eq:SD-num-1} and \( \chi_h := \lambda_h^m \) in \eqref{eq:SD-num-2}, and subsequently summing up equations \eqref{eq:SD-num-1}-\eqref{eq:SD-num-3}, we obtain
\begin{align}\label{thm:area-decreasing-mass-SD-2}
\int_{\Gamma_h^{m-1}} \nabla_{\Gamma_h^{m-1}} (\tau v_h^m + {\rm id}) \cdot \nabla_{\Gamma_h^{m-1}} v_h^m + \|\nabla_{\Gamma_h^{m-1}} \lambda_h^m\|_{L^2(\Gamma_h^{m-1})}^2 +\frac{1}{\alpha} | c^m|^2\|T_h^{m-1}\|_{L^2_h} = 0.
\end{align}
By substituting inequality \eqref{area-decreasing-important-property} into \eqref{thm:area-decreasing-mass-SD-2}, we obtain 
\begin{align}\label{thm:area-decreasing-mass-SD-3}
        |\Gamma_h^m| - |\Gamma_h^{m-1}| \le \int_{\Gamma_h^{m-1}} \nabla_{\Gamma_h^{m-1}} (\tau v_h^m + {\rm id}) \cdot \nabla_{\Gamma_h^{m-1}} v_h^m \le 0.
\end{align}
This completes the proof.
\hfill\end{proof}

The advantages of the BGN-MDR method for surface diffusion in improving the mesh quality are shown in the following example.

\begin{example}\label{Example2}\upshape
We consider surface diffusion with an initial surface being a box centered at \((0,0,0)\) with dimensions \(1 \times 6 \times 1\) in the \(x\)-, \(y\)-, and \(z\)-directions, respectively. The numerical results obtained by the BGN and BGN-MDR methods with mesh size \(h = 0.2\) and time step size \(\tau = 10^{-2}\) are shown in Figures \ref{Example2d} and \ref{Example2e}, where both schemes preserve good mesh quality. 

However, when the smaller time step size \(\tau = 10^{-4}\) is employed, the BGN scheme produces mesh distortion, whereas the BGN-MDR scheme maintains good mesh quality; see Figures~\ref{Example2b} and \ref{Example2c}. Furthermore, Figure~\ref{eg:Example2-auxiliary} shows the evolution of the surface area and \(\|T_h^{m-1}\|_{L^2_h}\). These numerical results demonstrate the superiority of the BGN-MDR method over the BGN method in improving mesh quality while maintaining energy stability. 


    
\begin{figure}[htbp]
        \centering
        \vspace{-5pt}
    
        \begin{subfigure}[b]{0.6\textwidth}
        \centering
            \includegraphics[width=0.7\textwidth]{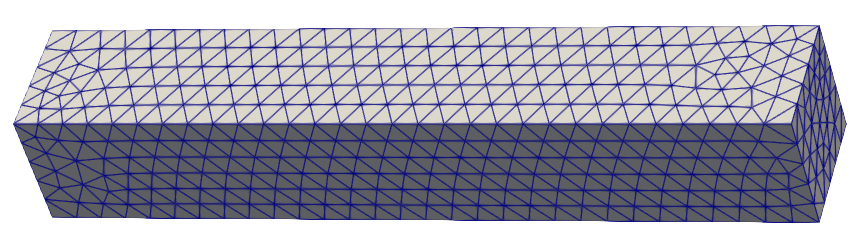}
            \caption{Initial Surface}
            \label{Example2a}
        \end{subfigure}
\vspace{8pt} 
    
        \begin{subfigure}[b]{0.45\textwidth}
        \centering
            \includegraphics[width=0.4\textwidth]{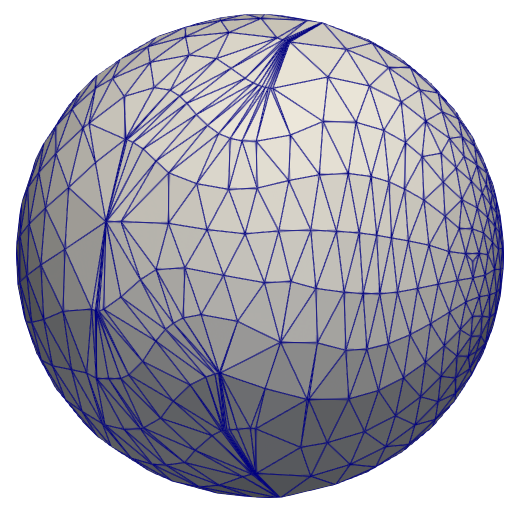}
            \caption{BGN scheme for $\tau=10^{-4}$}
            \label{Example2b}
        \end{subfigure}
        \hspace{10pt}
        \begin{subfigure}[b]{0.45\textwidth}
        \centering
            \includegraphics[width=0.4\textwidth]{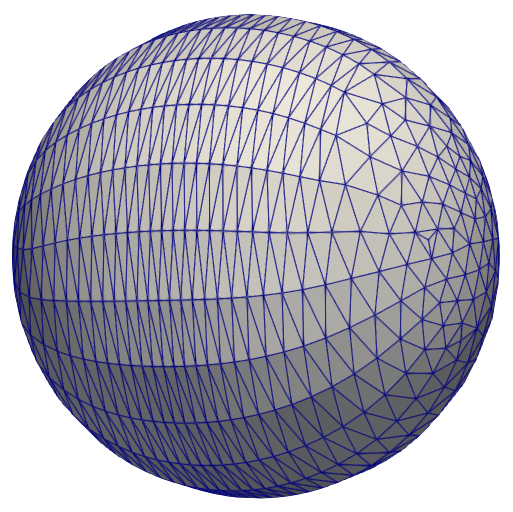}
            \caption{BGN-MDR scheme for $\tau=10^{-4}$}
            \label{Example2c}
        \end{subfigure}
    
        \vspace{8pt} 
    
        \begin{subfigure}[b]{0.45\textwidth}
        \centering
            \includegraphics[width=0.4\textwidth]{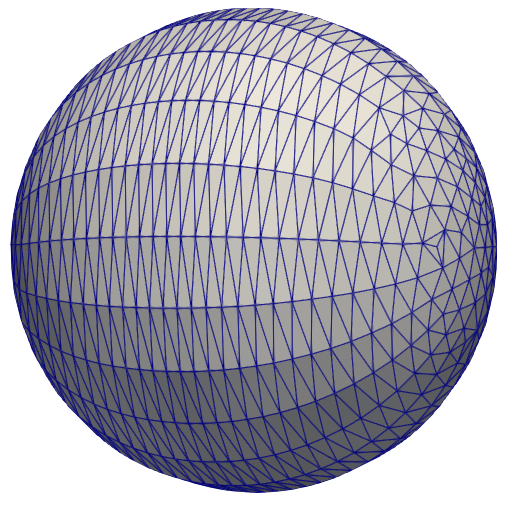}
            \caption{BGN scheme for $\tau=10^{-2}$}
            \label{Example2d}
        \end{subfigure}
        \hspace{10pt}
        \begin{subfigure}[b]{0.45\textwidth}
        \centering
            \includegraphics[width=0.4\textwidth]{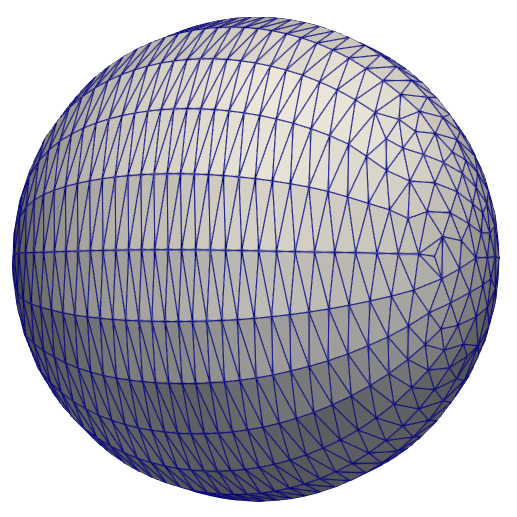}
            \caption{BGN-MDR scheme for $\tau=10^{-2}$}
            \label{Example2e}
        \end{subfigure}
    
        \vspace{-5pt}
        \caption{Surface evolution in Example \ref{Example2}}
        \label{eg:Example2}
\end{figure}

    \begin{figure}[htbp]
        \centering
%
        \begin{subfigure}[b]{0.45\textwidth}
            \includegraphics[width=0.9\textwidth]{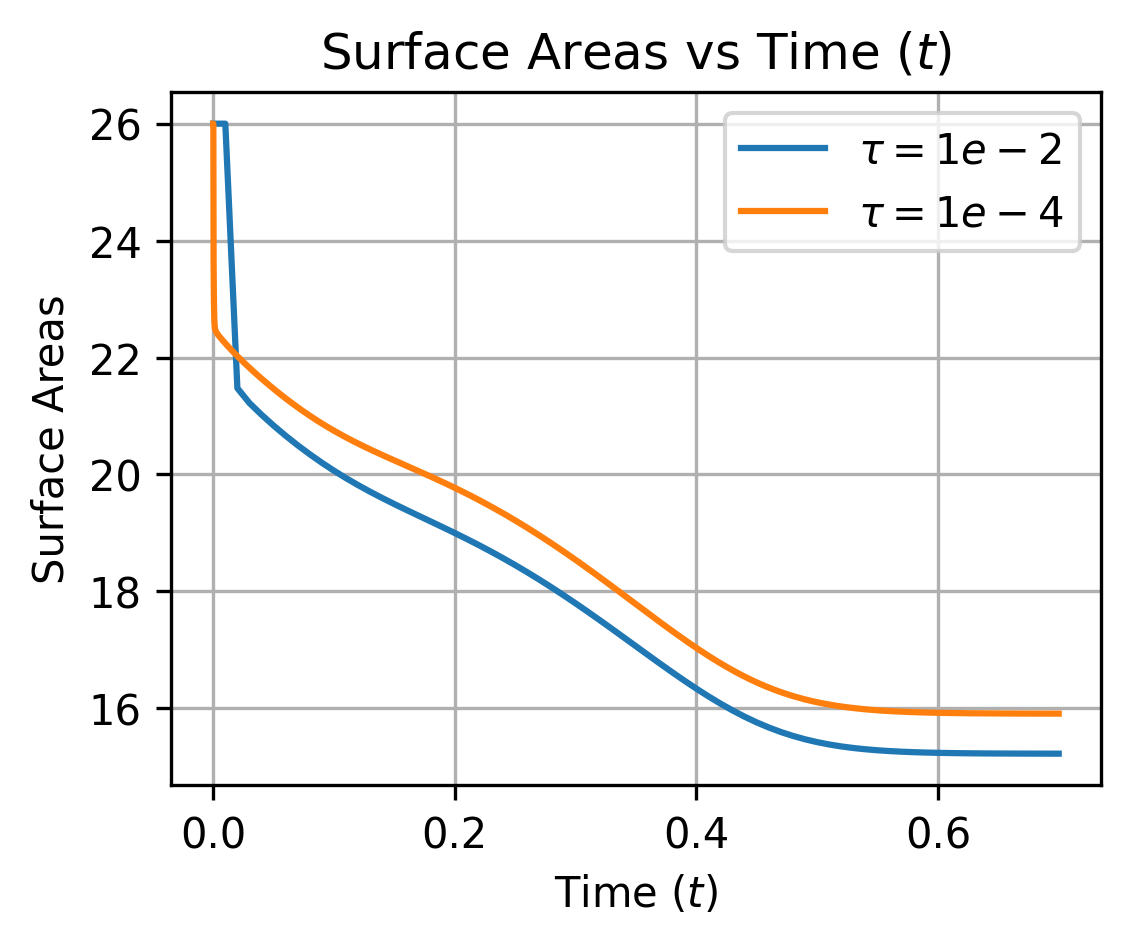}
            \caption{Surface areas versus time $t$}
            \label{Example2-area-1e-4}
        \end{subfigure}
        \hfill
%
        \begin{subfigure}[b]{0.45\textwidth}
            \includegraphics[width=0.9\textwidth]{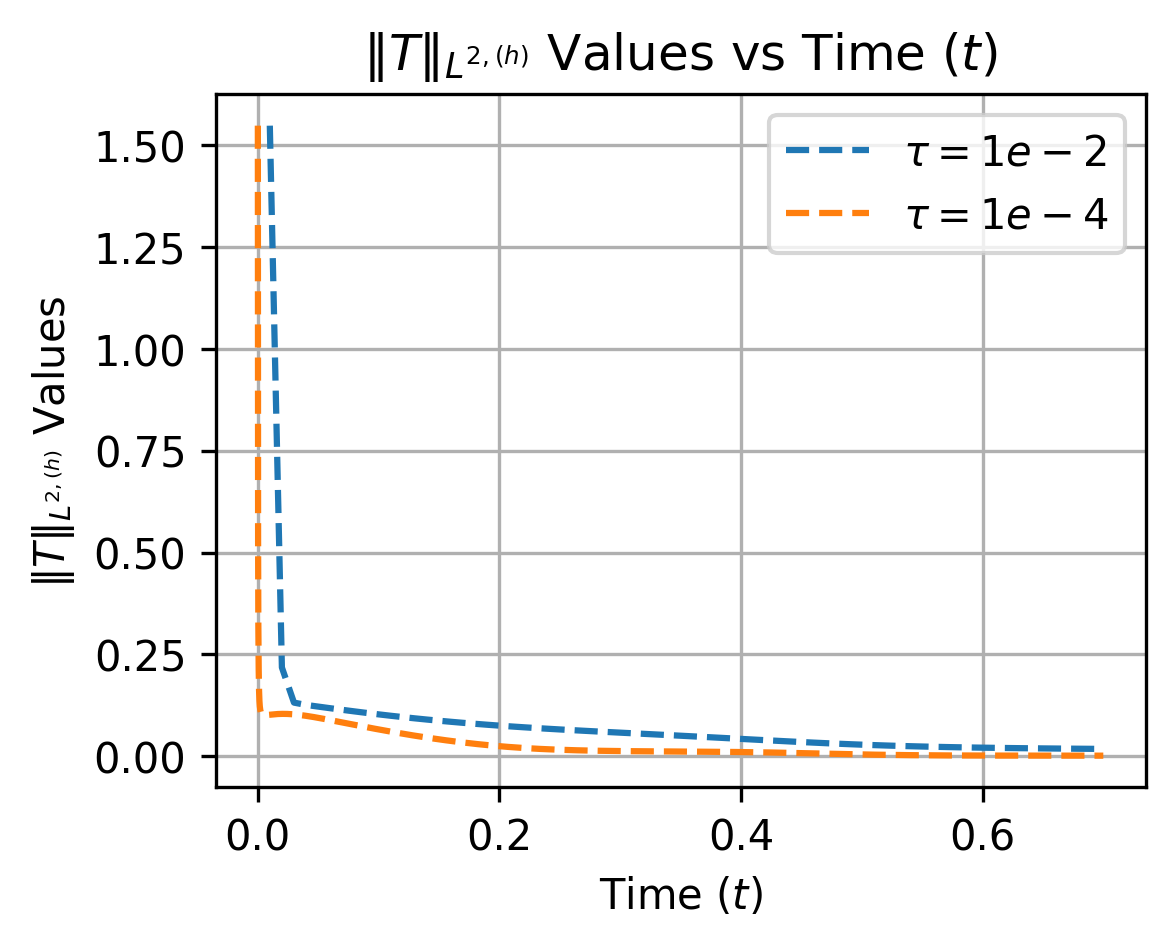}
            \caption{$\|T\|_{L^2_h}$ values versus time $t$}
            \label{Example2-T2-1e-4}
        \end{subfigure}
%
        \caption{Evolution of surface area and $\|T_h^{m-1}\|_{L^2_h}$ in Example \ref{Example1}.}
        \label{eg:Example2-auxiliary}
    \end{figure}
    
\end{example}

\begin{figure}[htbp]
        \centering
%
        \begin{subfigure}[b]{0.45\textwidth}
            \includegraphics[width=0.9\textwidth]{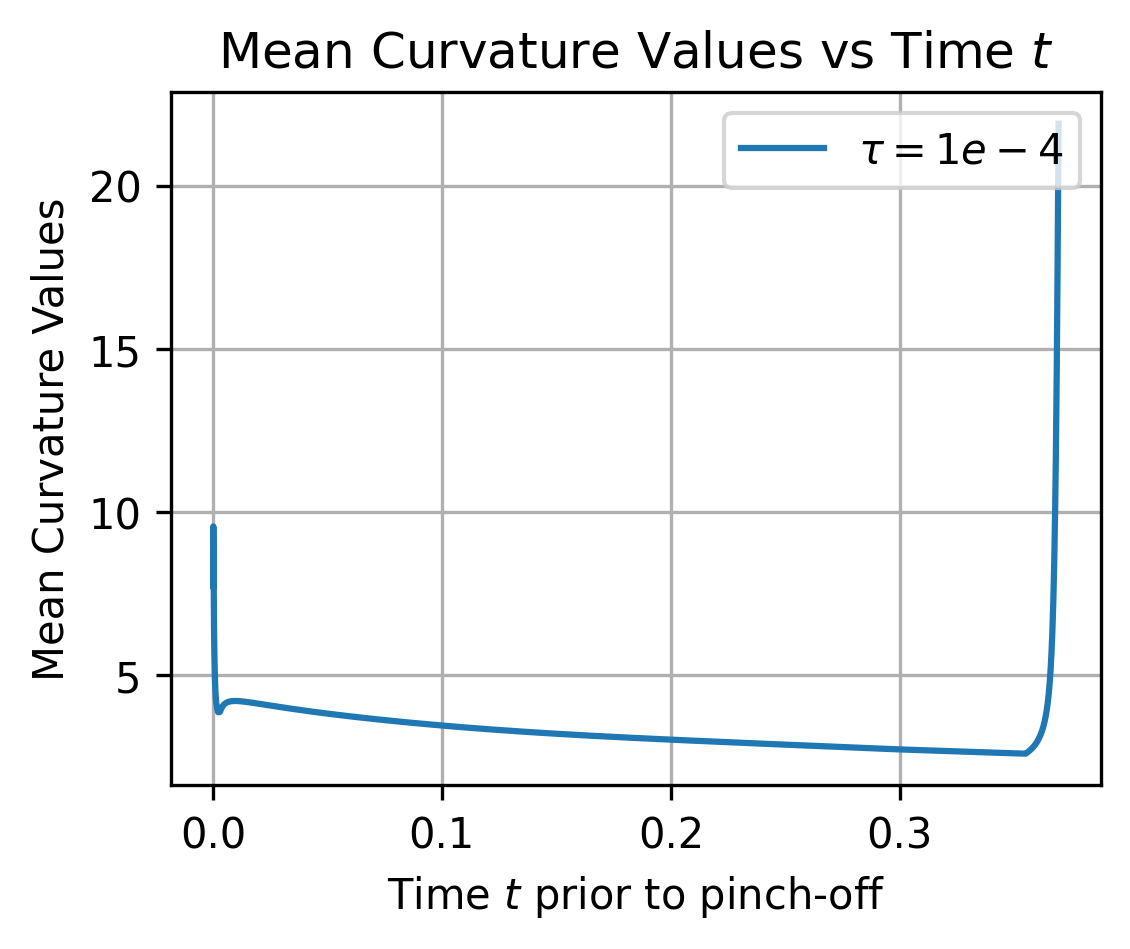}
            \caption{Maximum mean curvature values versus time \(t\) (prior to pinch-off).}
            \label{Example8-prior}
        \end{subfigure}
        \hfill
%
        \begin{subfigure}[b]{0.45\textwidth}
            \includegraphics[width=0.9\textwidth]{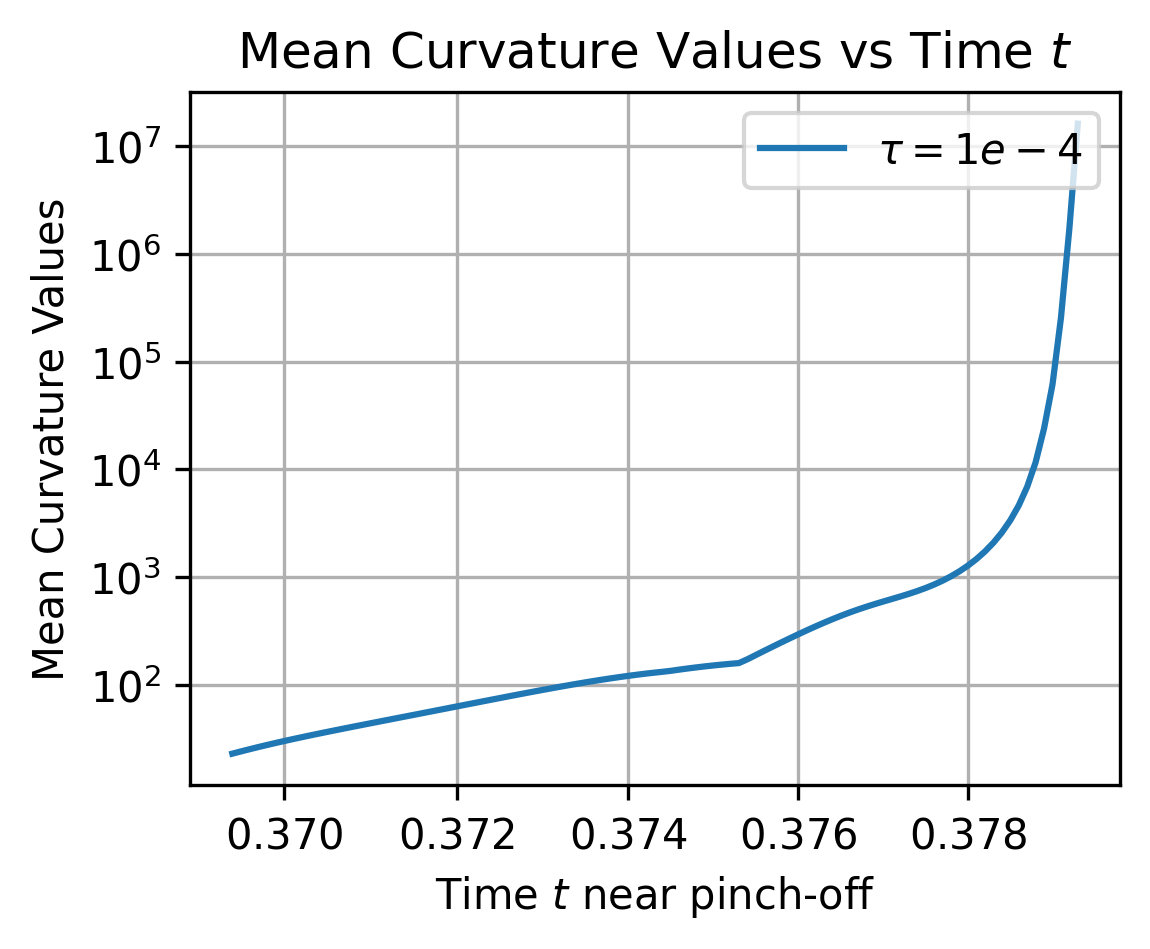}
            \caption{Maximum mean curvature values versus time $t$ (near pinch-off)}
            \label{Example8-near}
        \end{subfigure}
%
        \caption{{Evolution of maximum mean curvatures in Example \ref{Example1-1-8}.}}
        \label{eg:Example8-auxiliary}
    \end{figure}

\begin{example}\label{Example1-1-8}\upshape
{We also consider the surface diffusion flow starting from an initial surface given by a box centered at \((0,0,0)\) with dimensions \(1 \times 8 \times 1\) in the \(x\)-, \(y\)-, and \(z\)-directions, respectively. Numerical results obtained using the BGN and BGN-MDR methods with mesh size \(h = 0.2\) and time step size \(\tau = 10^{-3}\) are shown in Figures~\ref{Example8d} and~\ref{Example8e}. Both schemes preserve good mesh quality and are capable of predicting the pinch-off time associated with singularity formation.

However, when a smaller time step size \(\tau = 10^{-4}\) is used, the BGN scheme leads to mesh distortion, whereas the BGN-MDR scheme continues to preserve good mesh quality and still predicts the pinch-off time; see Figures~\ref{Example8b} and~\ref{Example8c}. Furthermore, Figure~\ref{eg:Example8-auxiliary} presents the evolution of the maximum mean curvature over the discrete surface and clearly demonstrates its blow-up behavior as the pinch-off time is approached during the formation of the singularity.

}

\end{example}

{ 
\begin{figure}[htbp]
        \centering
        \vspace{-5pt}
    
        \begin{subfigure}[b]{0.6\textwidth}
        \centering
            \includegraphics[width=1\textwidth]{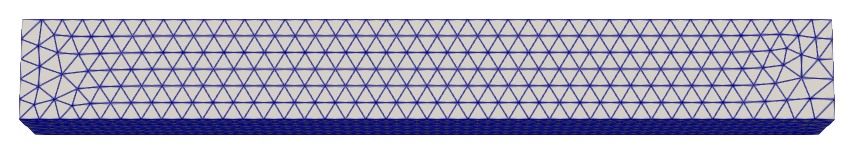}
            \caption{Initial Surface}
            \label{Example8a}
        \end{subfigure}
\vspace{8pt} 
    
        \begin{subfigure}[b]{0.45\textwidth}
        \centering
            \includegraphics[width=0.7\textwidth]{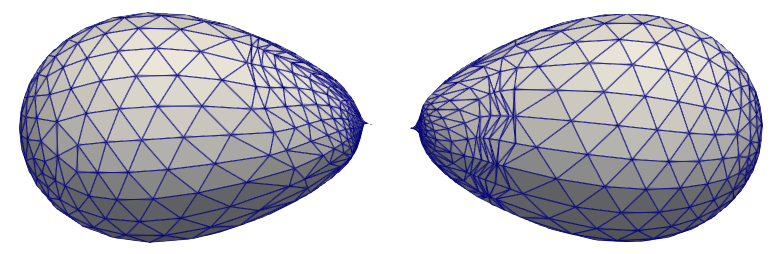}
            \caption{BGN scheme with $\tau=10^{-4}$\\
            \centering ($t = 0.3509$)}
            \label{Example8b}
        \end{subfigure}
        \hspace{10pt}
        \begin{subfigure}[b]{0.45\textwidth}
        \centering
            \includegraphics[width=0.7\textwidth]{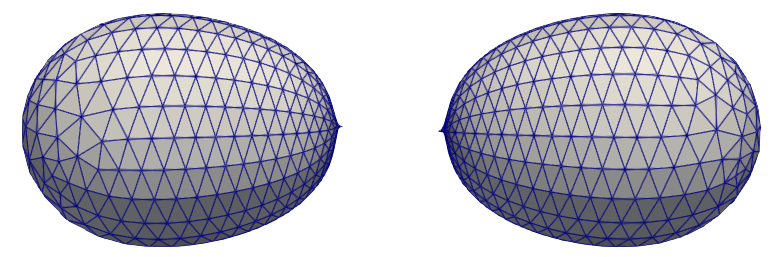}
            \caption{BGN-MDR scheme with $\tau=10^{-4}$\\
            \centering ($t = 0.3793$)}
            \label{Example8c}
        \end{subfigure}
    
        \vspace{8pt} 
    
        \begin{subfigure}[b]{0.45\textwidth}
        \centering
            \includegraphics[width=0.7\textwidth]{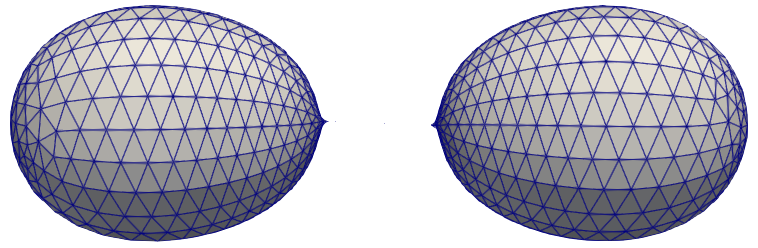}
            \caption{BGN scheme with $\tau=10^{-3}$\\
            \centering ($t = 0.369$)}
            \label{Example8d}
        \end{subfigure}
        \hspace{10pt}
        \begin{subfigure}[b]{0.45\textwidth}
        \centering
            \includegraphics[width=0.7\textwidth]{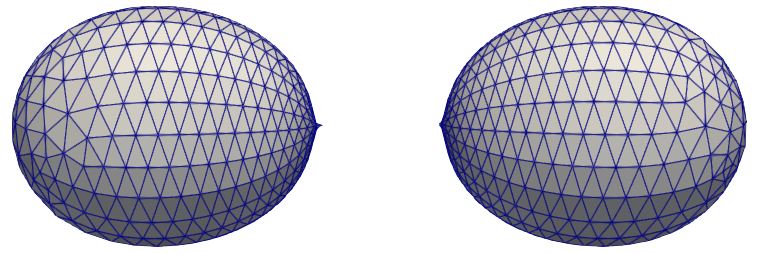}
            \caption{BGN-MDR scheme with $\tau=10^{-3}$\\
            \centering ($t = 0.389$)}
            \label{Example8e}
        \end{subfigure}
    
        \vspace{-5pt}
        \caption{{Surface evolution in Example \ref{Example1-1-8}}}
        \label{eg:Example8}
\end{figure}}
    
\section{BGN-MDR scheme for open surfaces with moving contact lines}\label{section:4}

\subsection{Introduction to solid-state dewetting (SSD)}
    

Solid-state dewetting (SSD) is an evolution phenomenon on open surfaces driven by surface diffusion and characterized by moving contact lines. As a thin film migrates along a substrate, a contact line forms at the junction of the solid, vapor, and substrate phases, adding extra kinetic complexity \cite{Thompson1,Armelao,Schmidt}.

In this section, we introduce the BGN-MDR scheme, which extends our previous approach for closed surfaces, for mean curvature flow and surface diffusion of open surfaces with moving contact lines.
    
\subsection{Mathematical model and governing equations in two dimensions (2D)}
Let the evolving curve
\[
    \Gamma(t) = X(s,t) = \begin{pmatrix} x(s,t) \\ y(s,t) \end{pmatrix},
\]
parameterized by the arc length \( s \in [0, L(t)] \), represent the film-vapor interface. Here, \( L(t) \) denotes the length of the curve, and \( x_{\rm c}^l(t) \) and \( x_{\rm c}^r(t) \) are the positions of the left and right moving contact points, respectively. The surface diffusion and mean curvature flow are governed by
\begin{align*}
&\mbox{Surface diffusion:}\quad 
&&\begin{cases}
\partial_t X \cdot n = \partial_s^2 \kappa, & 0 < s < L(t),\quad t \ge 0, \\
\kappa = -\bigl(\partial_s^2 X\bigr) \cdot n, & \quad n = \bigl(-\partial_s y,\, \partial_s x\bigr)^T,
\end{cases}\\[10pt] 
&\mbox{Mean curvature flow:}\quad 
&&\begin{cases}
\partial_t X \cdot n = -\kappa, & 0 < s < L(t),\quad t \ge 0, \\
\kappa = -\bigl(\partial_s^2 X\bigr) \cdot n, & \quad n = \bigl(-\partial_s y,\, \partial_s x\bigr)^T , 
\end{cases}
\end{align*}
where \( \kappa \) is the curvature and \( n \) is the unit normal vector. The initial condition is
    \[
    X(s,0) = X_0(s) = \begin{pmatrix} x_0(s) \\ y_0(s) \end{pmatrix}, \quad s \in [0, L_0].
    \]
The boundary conditions are as follows. At the contact points, the interface satisfies the contact point condition
    \[
    y(0,t) = 0,\quad y(L,t) = 0,\quad t \ge 0,
    \]
    and Young's law (contact angle condition)
    \[
    \cos\theta_{\rm d}^l = \sigma = \cos\theta_{\rm d}^r,\quad t \ge 0,
    \]
    where \(\theta_{\rm d}^l\) and \(\theta_{\rm d}^r\) are the angles between the interface and the substrate at the left and right contact points, respectively. For surface diffusion, an additional zero-flux condition is imposed:
    \[
    \partial_s \kappa(0,t) = 0,\quad \partial_s \kappa(L,t) = 0,\quad t \ge 0.
    \]

For any test function 
\(
\eta := (\eta_1, \eta_2) \in H^1(\Gamma(t)) \times H^1(\Gamma(t)),
\)
the following relation holds from integration by parts:
\begin{align}\label{eq:ibp}
    (\kappa n, \eta)_{\Gamma(t)} &= (-\partial_s^2 X, \eta)_{\Gamma(t)} = (\partial_s X, \partial_s \eta)_{\Gamma(t)} - (\partial_s X \cdot \eta)\Big|_{s=0}^{s=L(t)} \\
    &= (\partial_s X, \partial_s \eta)_{\Gamma(t)} - \sigma\,[\eta_1^r - \eta_1^l] - \sqrt{1-\sigma^2}\,[\eta_2^r - \eta_2^l],
\end{align}
where the last equality is obtained by applying Young's law and the pairs \((\eta_1^r,\eta_2^r)\) and \((\eta_1^l,\eta_2^l)\) denote the values of \(\eta\) at the right and left moving contact points, respectively.

Correspondingly, a functional \( L_h^{m-1}(\eta_h) \) is defined on the finite element function space \((S_h^{m-1})^2\) under spatial discretization by
\[
L_h^{m-1}(\eta_h) := \bigl( \partial_s X_h^{m-1}, \partial_s \eta_h \bigr) - \sigma\,[\eta_{h,1}^{r,m-1} - \eta_{h,1}^{l,m-1}]-\sqrt{1-\sigma^2} \,[\eta_{h,2}^{r,m-1} - \eta_{h,2}^{l,m-1}].
\]
By virtue of the Riesz representation theorem, it is deduced that there exists \(\nu_h^{m-1} \in (S_h^{m-1} )^2\) such that
\[
L_h^{m-1}(\eta_h) = \bigl( \nu_h^{m-1}, \eta_h \bigr)^{(h)} \quad \text{for all } \eta_h \in (S_h^{m-1})^2.
\]

Then the tangential vector \( T_h^{m-1} \) is defined similarly to the equation \eqref{eq:T_h^m}, and it is required to determine \((T_h^{m-1}, \lambda_h) \in (S_h^{m-1})^2 \times S_h^{m-1}\) such that
\begin{equation}\label{eq:T_h^m-higher-order-2D}
\begin{aligned}
\bigl( T_h^{m-1}, \eta_h \bigr)^{(h)} + \bigl( \lambda_h\, n_h^{m-1}, \eta_h \bigr)^{(h)} &= \bigl( \nu_h^{m-1}, \eta_h \bigr)^{(h)}, \quad \forall \eta_h \in (S_h^{m-1})^2, \\[1mm]
\bigl( T_h^{m-1}, n_h^{m-1}\, \chi_h \bigr)^{(h)} &= 0, \quad \forall \chi_h \in S_h^{m-1}.
\end{aligned}
\end{equation}
In this formulation, \( T_h^{m-1} \in (S_h^{m-1})^2 \) is understood to represent the component of \(\nu_h^{m-1}\) that is orthogonal to the space \(\mathrm{span}\{ I_h(n_h^{m-1}\, \chi_h) \mid \chi_h \in S_h^{m-1} \}\).

\subsection{Numerical scheme for mean curvature flow and surface diffusion in 2D}\label{sec:num-SSD-2D}
    Let \( \mathring{S}_h(\Gamma_h^{m-1}) \) denote the subspace of \( S_h^{m-1} \) comprising finite element functions that vanish on \( \partial\Gamma_h^{m-1} \), and define \( \mathbf{X}_h^{m-1} := S_h^{m-1} \times \mathring{S}_h(\Gamma_h^{m-1}) \); the discrete velocity \( v_h^m \) of the flow map \( X_h^m \) is accordingly assumed to lie in \( \mathbf{X}_h^{m-1} \) so as to ensure that the contact point condition is satisfied. By emulating the BGN-MDR scheme \eqref{eq:MCF-num-equ} for closed surfaces, the following numerical scheme for mean curvature flow of open curves with moving contact points is proposed: Find \( (v_h^m, \lambda_h^m,  c^m) \in \mathbf{X}_h^{m-1} \times S_h^{m-1} \times \mathbb{R} \) such that
    \begin{subequations}\label{eq:MCF-num-2D}
        \begin{align}
            (\partial_s (\tau v_h^m +{\rm id}), \partial_s \eta_h) - \bigl( \lambda_h^mn_h^{m-1}, \eta_h\bigr)^{(h)} - c^m \bigl( T_h^{m-1}, \eta_h\bigr)^{(h)} &= \sigma \bigl[\eta_{h,1}^{r,m-1} - \eta_{h,1}^{l,m-1}\bigr], \label{eq:MCF-num-1-2D} \\[1mm]
            \bigl( v_h^m \cdot {n}_h^{m-1}, \chi_h\bigr)^{(h)} + \bigl( \lambda_h^m, \chi_h\bigr)^{(h)} &= 0,\label{eq:MCF-num-2-2D} \\[1mm]
            \bigl( v_h^m, T_h^{m-1}\bigr)^{(h)} + \alpha c^m \|T_h^{m-1}\|_{L^2_h} &= 0 \label{eq:MCF-num-3-2D}
        \end{align}    
    \end{subequations}
for all \( (\eta_h, \chi_h) \in \mathbf{X}_h^{m-1} \times S_h^{m-1} \). By solving for \(\lambda_h^m\) and \(c^m\) from \eqref{eq:MCF-num-2-2D} and \eqref{eq:MCF-num-3-2D}, the above scheme \eqref{eq:MCF-num-2D} may be equivalently reformulated in a manner that closely parallels the numerical scheme in \eqref{eq:MCF-scheme--equal}:
    \begin{subequations}\label{eq:MCF-num-2D-proj}
        \begin{align}
            (\partial_s (\tau v_h^m + {\rm id}), \partial_s \eta_h) &+ \bigl( v_h^m \cdot n_h^{m-1}, \eta_h\cdot n_h^{m-1}\bigr)^{(h)} + \frac{\bigl( T_h^{m-1}, \eta_h\bigr)^{(h)} \bigl( T_h^{m-1}, v_h^m\bigr)^{(h)}}{\alpha\|T_h^{m-1}\|_{L^2_h}} \notag\\[1mm]
            &= \sigma \bigl[\eta_{h,1}^{r,m-1} - \eta_{h,1}^{l,m-1}\bigr], \quad \text{for all } \eta_h \in \mathbf{X}_h^{m-1}.
        \end{align}
    \end{subequations}
    
    Furthermore, the BGN-MDR numerical scheme for surface diffusion of open curves with moving contact points is formulated as follows: Find \( (v_h^m, \lambda_h^m, c^m) \in \mathbf{X}_h^{m-1} \times S_h^{m-1} \times \mathbb{R} \) such that
    \begin{subequations}\label{eq:SD-num-2D}
        \begin{align}
            (\partial_s (\tau v_h^m + {\rm id}), \partial_s \eta_h) - \bigl( \lambda_h^m n_h^{m-1}, \eta_h\bigr)^{(h)} -  c^m \bigl( T_h^{m-1}, \eta_h\bigr)^{(h)} &= \sigma \bigl[\eta_{h,1}^{r,m-1} - \eta_{h,1}^{l,m-1}\bigr],  \label{eq:SD-num-1-2D} \\[1mm]
            \bigl( v_h^m \cdot n_h^{m-1}, \chi_h\bigr)^{(h)} + \bigl( \partial_s \lambda_h^m, \partial_s \chi_h\bigr) &= 0, \label{eq:SD-num-2-2D} \\[1mm]
            \bigl( v_h^m, T_h^{m-1}\bigr)^{(h)} + \alpha c^m \|T_h^{m-1}\|_{L^2_h} &= 0  \label{eq:SD-num-3-2D}
        \end{align}    
    \end{subequations}
for all \( (\eta_h, \chi_h) \in \mathbf{X}_h^{m-1} \times S_h^{m-1} \). In an analogous fashion to the proofs of well-posedness for the BGN-MDR numerical schemes for closed surfaces, as established in Theorems \ref{thm:well-pose-mass-mcf} and \ref{thm:SD-num}, the well-posedness of the schemes for open curves with moving contact points can similarly be rigorously verified under the same mild conditions, since the boundary term 
    \[
    \sigma \bigl[\eta_{h,1}^{r,m-1} - \eta_{h,1}^{l,m-1}\bigr]
    \]
    in \eqref{eq:MCF-num-1-2D} or \eqref{eq:SD-num-1-2D} is eliminated in the corresponding homogeneous system. 
    
    \begin{theorem}\label{thm:SD-num-2D}
        Assume that the discrete curve $\Gamma_h^{m-1}$ satisfies the following mild conditions:
    \begin{enumerate}[label={(A\arabic*)}]
        \item The elements are nondegenerate, i.e., for each $K\subset \mathcal{K}_h^{m-1}$, it holds
      $
        |K| > 0.
       $
        \item The normal vectors 
        fulfill
        \[
        \dim \bigl(\operatorname{span}\big\{(n_h^{m-1}\phi_h,1)^{(h)}\big|\,\phi_h \in S_h^{m-1}\big\} \bigr) \;=\; 2.
        \]
        \end{enumerate}
        
        \noindent
        Then the BGN-MDR schemes in \eqref{eq:MCF-num-2D} and \eqref{eq:SD-num-2D} are well-posed, guaranteeing the existence of a unique solution \( (v_h^m, \lambda_h^m, c^m) \in \mathbf{X}_h^{m-1} \times S_h^{m-1} \times \mathbb{R} \).
    \end{theorem}
    
Moreover, the BGN-MDR numerical scheme maintains energy stability for both mean curvature flow and surface diffusion. More precisely, the following theorem holds.
    
    \begin{theorem}\label{thm:area-decreasing-mass-sd-2D}  
        Let \( (v_h^m, \lambda_h^m,  c^m) \in \mathbf{X}_h^{m-1} \times S_h^{m-1} \times \mathbb{R} \) be the solution of \eqref{eq:SD-num-2D}; then, the discrete energy 
        \[
        W_h^m := |\Gamma_h^m| - \sigma\bigl(x_r^m - x_l^m\bigr)
        \]
        is non-increasing, in the sense that
        \[
        W_h^m \le W_h^{m-1} \quad \text{for } m = 1,2,\ldots,N.
        \]
Here, \(x_l^m\) and \(x_r^m\) denote, respectively, the \(x\)-axis coordinates of the left-moving and right-moving contact points at time level \(t^m\). 
    \end{theorem}
    
    \begin{proof}
        Choosing test functions \( \eta_h := v_h^m \) in \eqref{eq:SD-num-1-2D} and \( \chi_h := \lambda_h^m \) in \eqref{eq:SD-num-2-2D} yields
\begin{align*}
        &(\partial_s (\tau v_h^m + {\rm id}), \partial_s v_h^m) + (\partial_s \lambda_h^m, \partial_s \lambda_h^m) + \frac{1}{\alpha}s(c^m)^2 \|T_h^{m-1}\|_{L^2_h} \\
        &= \frac{\sigma}{\tau} \Bigl[(x_c^r(t_m)-x_c^l(t_m)) - (x_c^r(t_{m-1})-x_c^l(t_{m-1}))\Bigr].
\end{align*}
        Then, by substituting inequality \eqref{area-decreasing-important-property} into the above relation, we obtain 
\[
        W_h^m - W_h^{m-1} \le (\partial_s X_h^m, \partial_s (X_h^m - X_h^{m-1})) - \sigma\Bigl[(x_c^r(t_m)-x_c^l(t_m)) - (x_c^r(t_{m-1})-x_c^l(t_{m-1}))\Bigr] \le 0,
\]
        which completes the proof.
  \hfill  \end{proof}
    
    The energy stability of the BGN-MDR numerical scheme for mean curvature flow is similarly established by an analogous argument. 
    
    The convergence and performance of the BGN-MDR method for mean curvature flow of open curves with moving contact points are tested in the following two examples. 
    
\begin{example}\label{example_para}\upshape
We consider the evolution of a curve in mean curvature flow with two endpoints constrained to the lines \(x = -\pi/4\) and \(x = \pi/4\), with the initial curve \(\Gamma^0\) parametrized by 
        \[
        (x(\theta), y(\theta)) = \Bigl(\theta,\; -\ln (\cos\theta) + 2 \Bigr), \quad \theta \in \left[-\frac{\pi}{4}, \frac{\pi}{4}\right].
        \]
The parametrization of the exact solution is given by $(x(\theta), y(\theta) +t)$, as shown in \cite{Garcke-phase}. 
{Since the constraints on the contact points are now specified on two parallel lines, rather than on the \(x\)-axis, our numerical scheme must be adapted accordingly. This involves updating the constraints in the function space to reflect the new locations of the contact points, as well as adjusting the boundary terms in the relation \eqref{eq:ibp} through integration by parts. Consequently, an analogous weak form for the derivation of the vector \(T_h^{m-1}\) and the corresponding numerical scheme for mean curvature flow can be constructed for the case where the contact points are located on \(x = \pm \pi/4\).
}

Figure~\ref{example_para_graph}(a) illustrates the evolution of the curve obtained with 31 mesh points and a time step size of \(\tau = 10^{-3}\). Moreover, Figure~\ref{example_para_graph}(b) shows that the distance errors between numerical and exact curves attains second-order convergence.

    \begin{figure}[htbp]
        \centering
        \begin{subfigure}[b]{0.45\textwidth}
            \includegraphics[width=\textwidth]{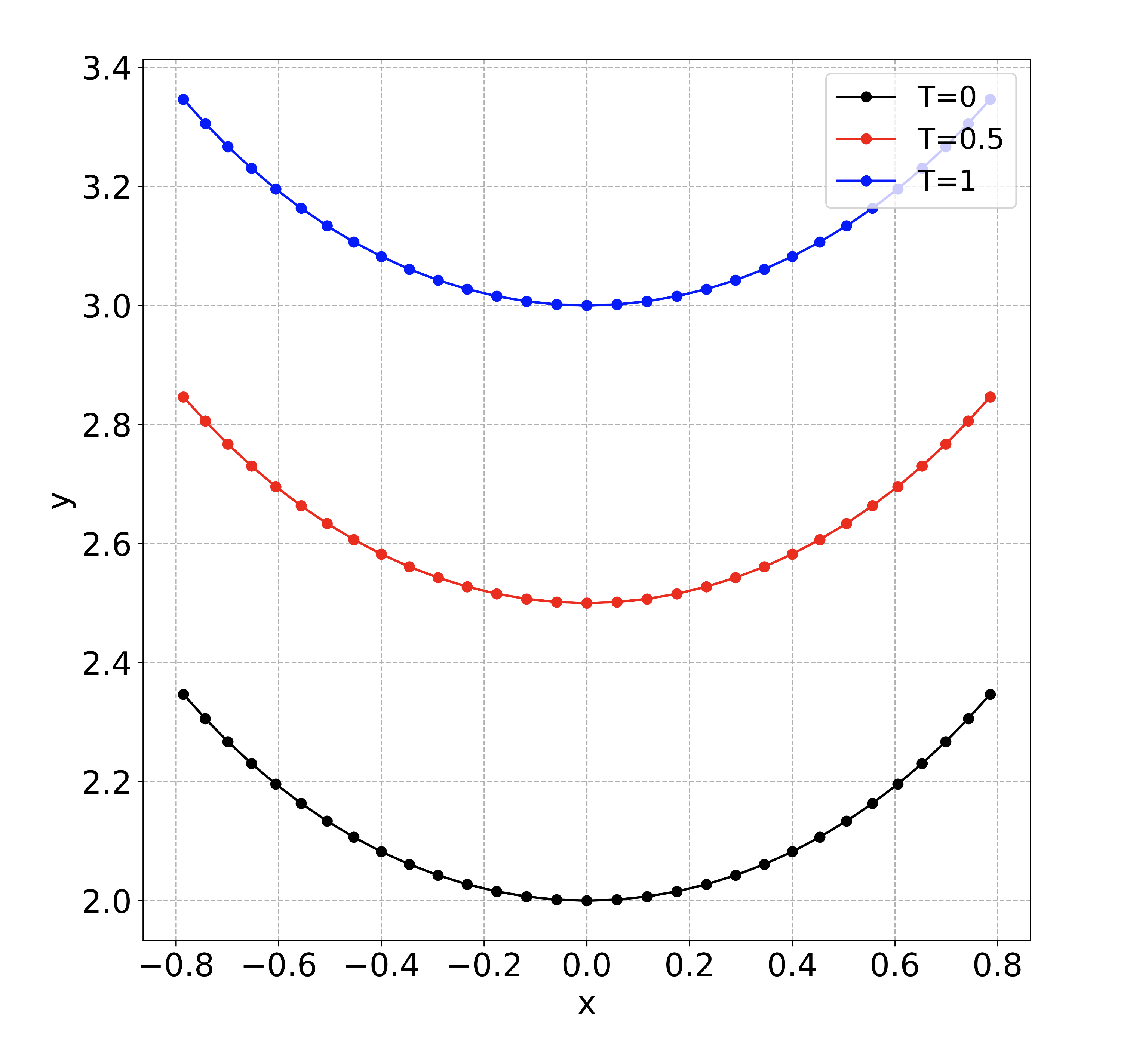}
            \caption{Evolution of the curve}
            \label{Example_para_curve1}
        \end{subfigure}
        \hfill
        \begin{subfigure}[b]{0.43\textwidth}
            \includegraphics[width=\textwidth]{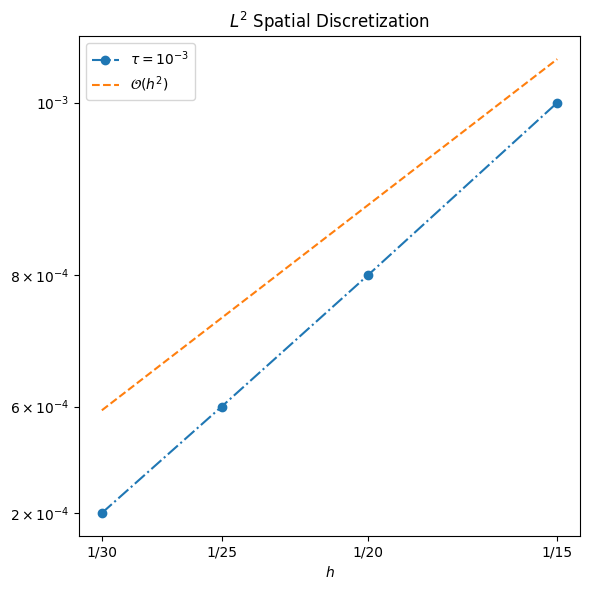}
            \caption{Errors from spatial discretization}
            \label{Example_para_curve2}
        \end{subfigure}
    
        \caption{Numerical results in Example~\ref{example_para}.}
        \label{example_para_graph}
    \end{figure}

\begin{example}\label{example_half_circle}\upshape
We consider the evolution of a curve in mean curvature flow with the two contact points constrained to the \(x\)-axis, with initial curve being a unit half circle. It is known that the exact evolving curve remains a half circle with radius
    $
    r(t)=\sqrt{1-2t}.
    $
The numerical solution given by the BGN-MDR scheme is presented in Figure \ref{fig:circle_rate_example}, which demonstrates that the mesh quality is well preserved during the evolution from \(T=0\) to \(T=0.2\).
    \begin{figure}[htbp]
        \centering
            \includegraphics[width=0.45\textwidth]{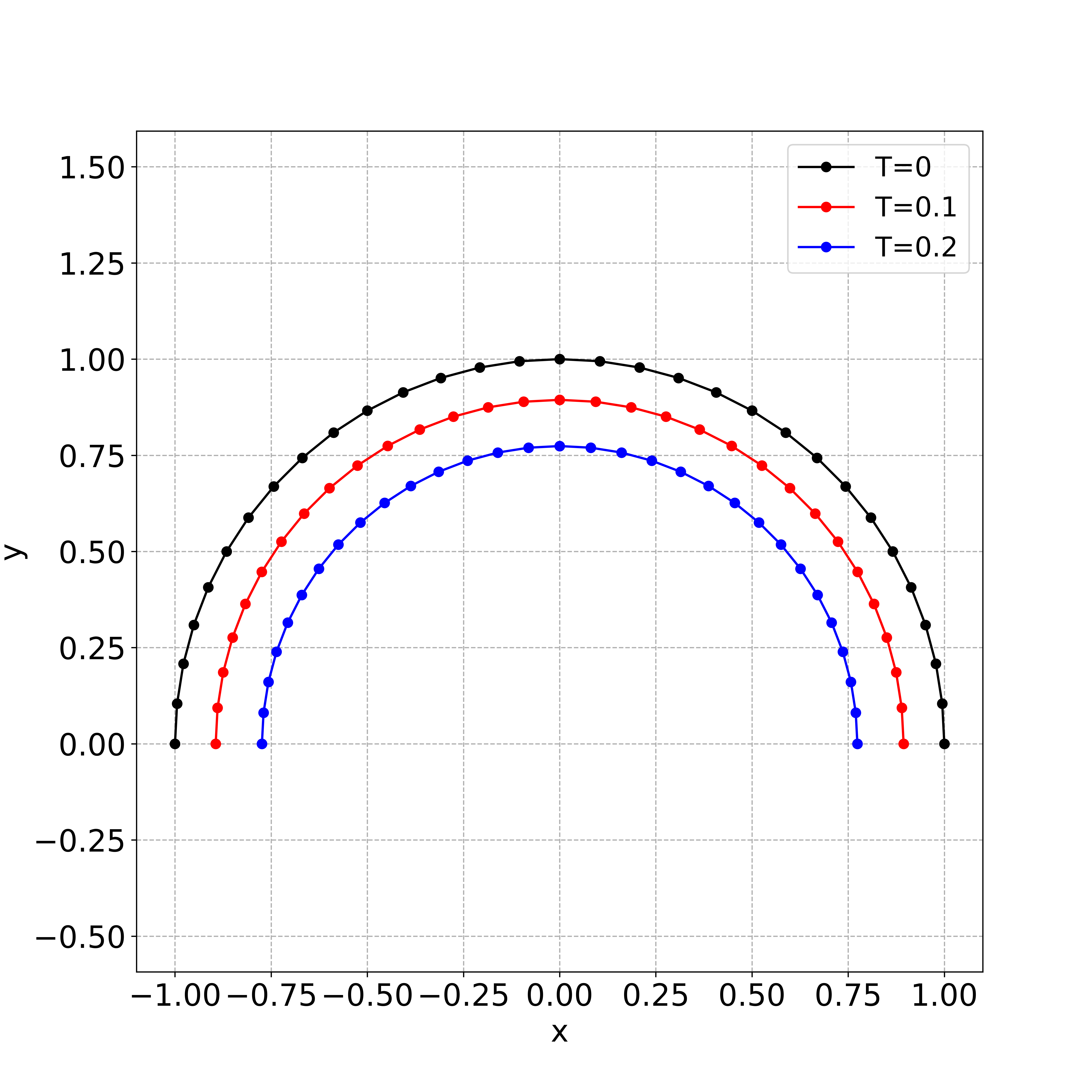}
        \caption{Evolution of the curve in Example~\ref{example_half_circle}.}
        \label{fig:circle_rate_example}
    \end{figure}
    
{Moreover, Figures~\ref{example_circle_rates} and~\ref{example_circle_rates_que6_H} present the distance errors and curvature errors from the spatial and temporal discretizations, respectively, demonstrating that the method attains the desired convergence rates in both space and time. Figure~\ref{example_circle_rates_que6_c} shows the auxiliary $c^m$ values for different mesh sizes. It can be observed that the $c^m$ values decrease as the mesh is refined, indicating that $c^m \to 0$ as \(h \to 0\). Figure~\ref{example_circle_rates_que6_evol} illustrates that, even when starting from a nonuniform mesh, the proposed numerical scheme leads to increasingly uniform mesh point distribution as time evolves.}

    \begin{figure}[htbp]
        \centering
        \begin{subfigure}[b]{0.45\textwidth}
            \includegraphics[width=\textwidth]{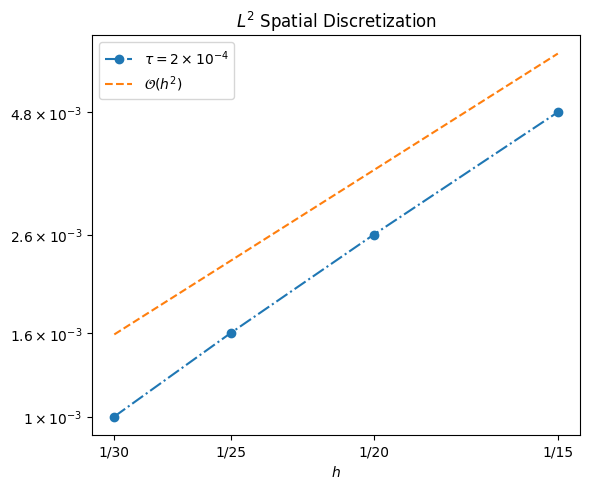}
            \caption{Spatial discretization distance errors}
            \label{Example_circle_spat}
        \end{subfigure}
        \hfill
        \begin{subfigure}[b]{0.45\textwidth}
            \includegraphics[width=\textwidth]{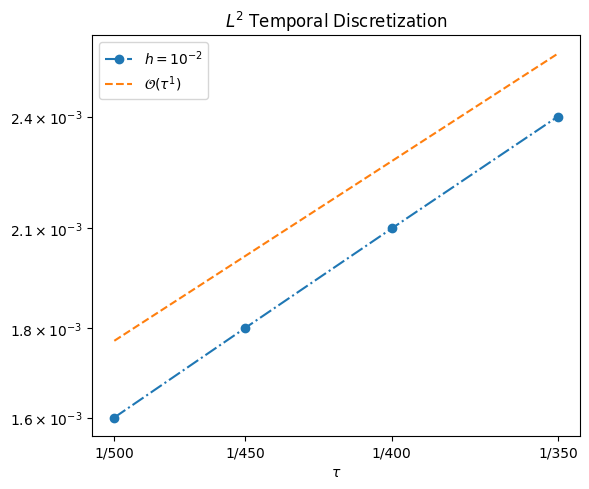}
            \caption{Temporal discretization distance errors}
            \label{Example_circle_temp}
        \end{subfigure}
    
        \caption{$L^\infty(0, 0.2; L^2)$ norm of the distance errors between the numerical and exact curves in Example~\ref{example_half_circle}.}
        \label{example_circle_rates}
    \end{figure}
    
    \begin{figure}[htbp]
        \centering
        \begin{subfigure}[b]{0.45\textwidth}
            \includegraphics[width=\textwidth]{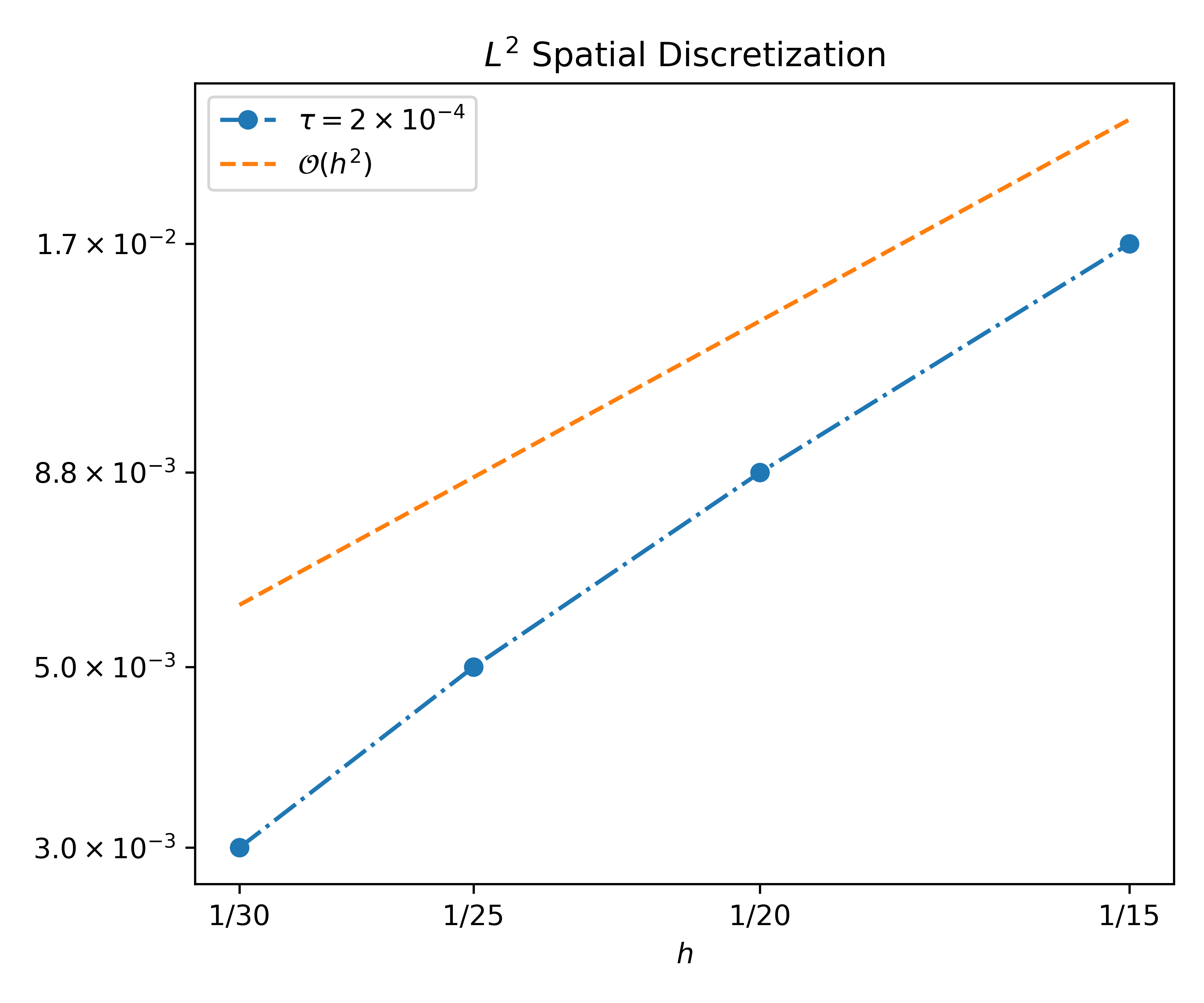}
            \caption{Spatial discretization curvature errors}
            \label{Example_circle_que6_H}
        \end{subfigure}
        \hfill
        \begin{subfigure}[b]{0.45\textwidth}
            \includegraphics[width=\textwidth]{./figs/circle_rate/L2-tem_half_circ.png}
            \caption{Temporal discretization curvature errors}
            \label{Example_circle_temp_que6_H}
        \end{subfigure}
    
        \caption{{\(L^\infty(0,0.2;L^2)\) norm of the curvature errors between the numerical and exact curves in Example~\ref{example_half_circle}.}}
        \label{example_circle_rates_que6_H}
    \end{figure}

      \begin{figure}[htbp]
        \centering
        \begin{subfigure}[b]{0.45\textwidth}
            \includegraphics[width=\textwidth]{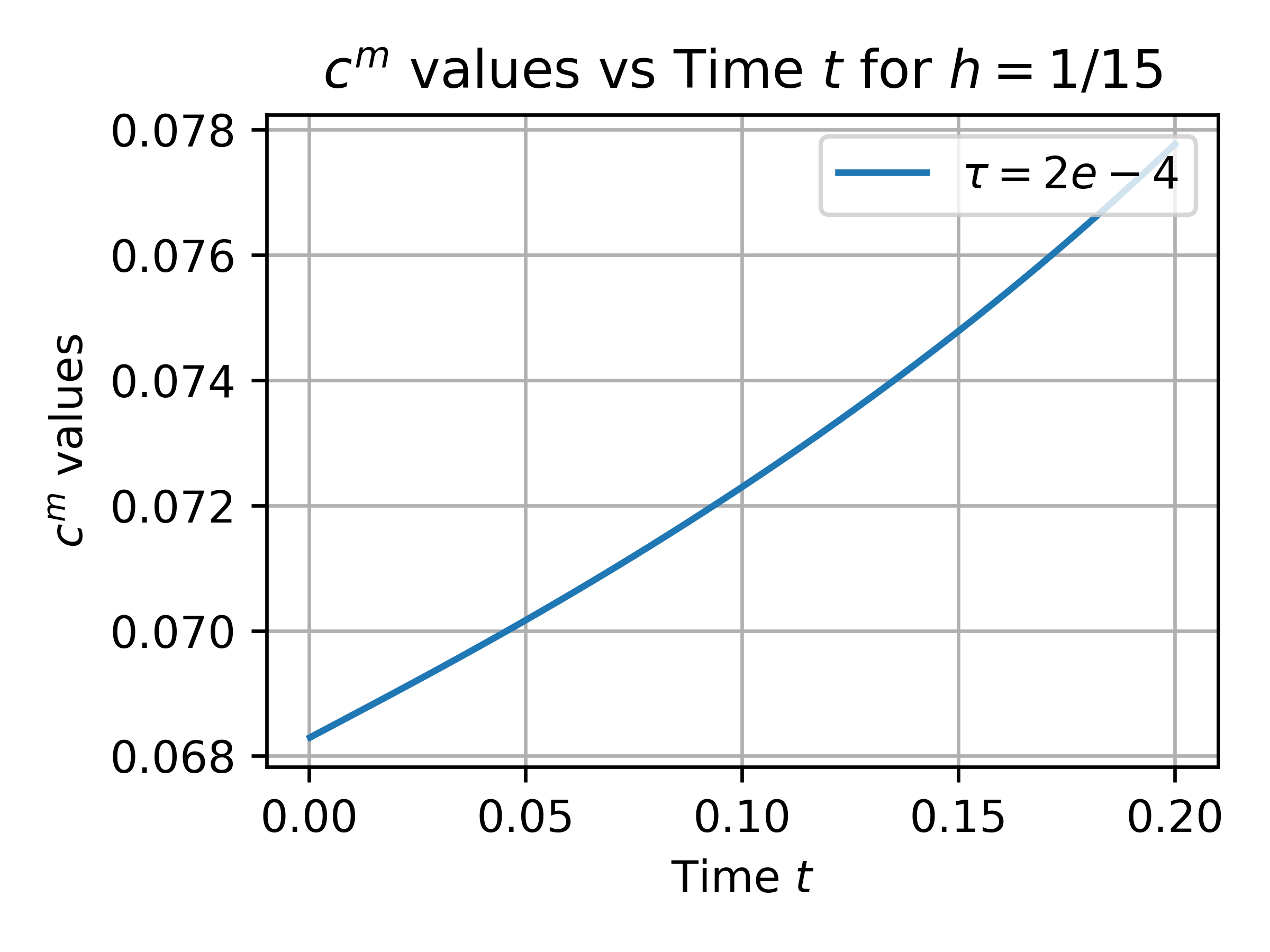}
            \caption{\(c^m\) versus time \(t\) (coarse mesh)}
            \label{Example_circle_que6_c}
        \end{subfigure}
        \hfill
        \begin{subfigure}[b]{0.45\textwidth}
            \includegraphics[width=\textwidth]{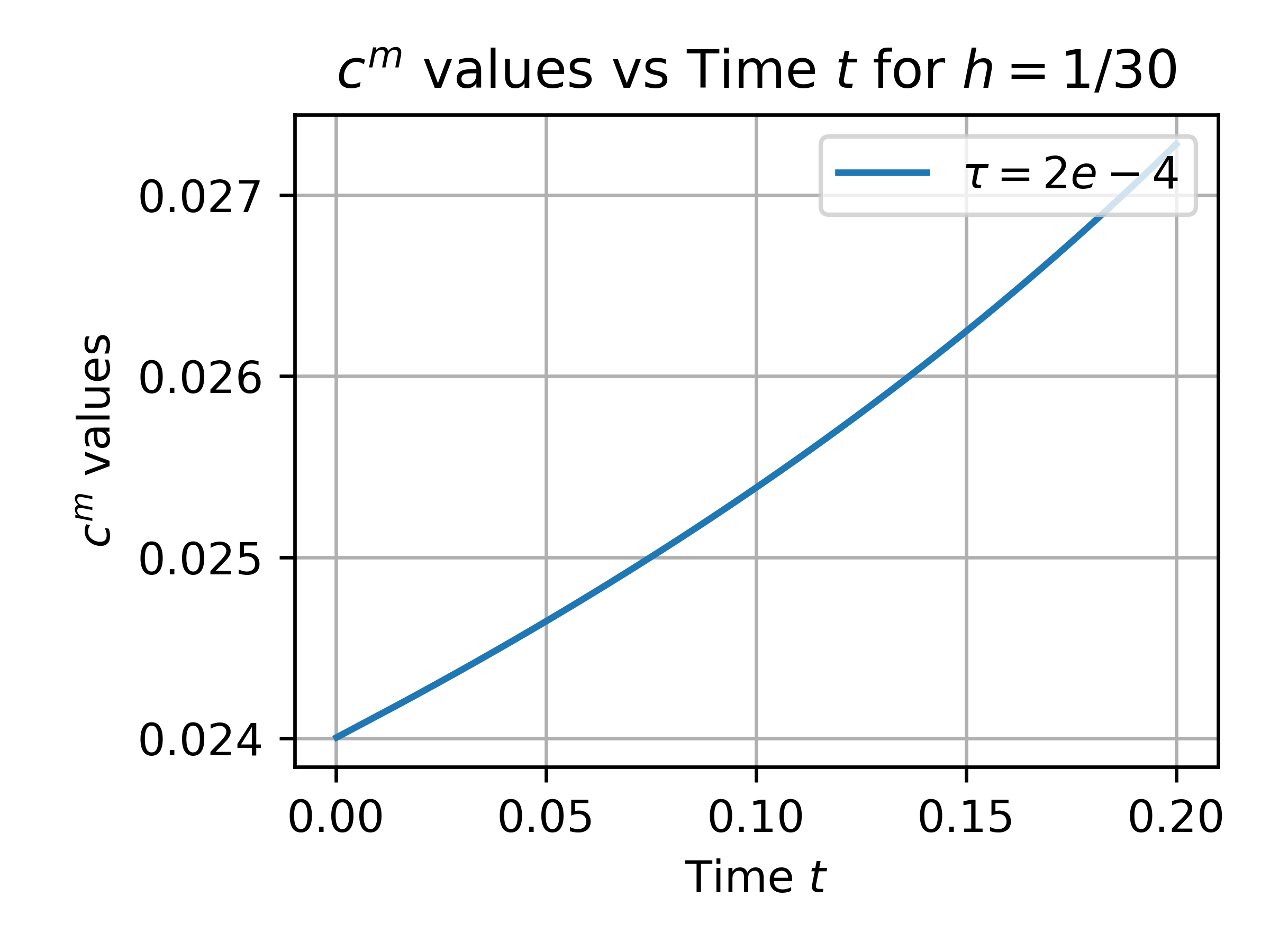}
            \caption{\(c^m\) versus time \(t\) (fine mesh)}
            \label{Example_circle_temp_que6_c}
        \end{subfigure}
    
        \caption{{Comparison of $c^m$ versus time $t$ for coarse and fine meshes in Example~\ref{example_half_circle}.}}
        \label{example_circle_rates_que6_c}
    \end{figure}

    \begin{figure}[htbp]
        \centering
        \begin{subfigure}[b]{0.45\textwidth}
            \includegraphics[width=\textwidth]{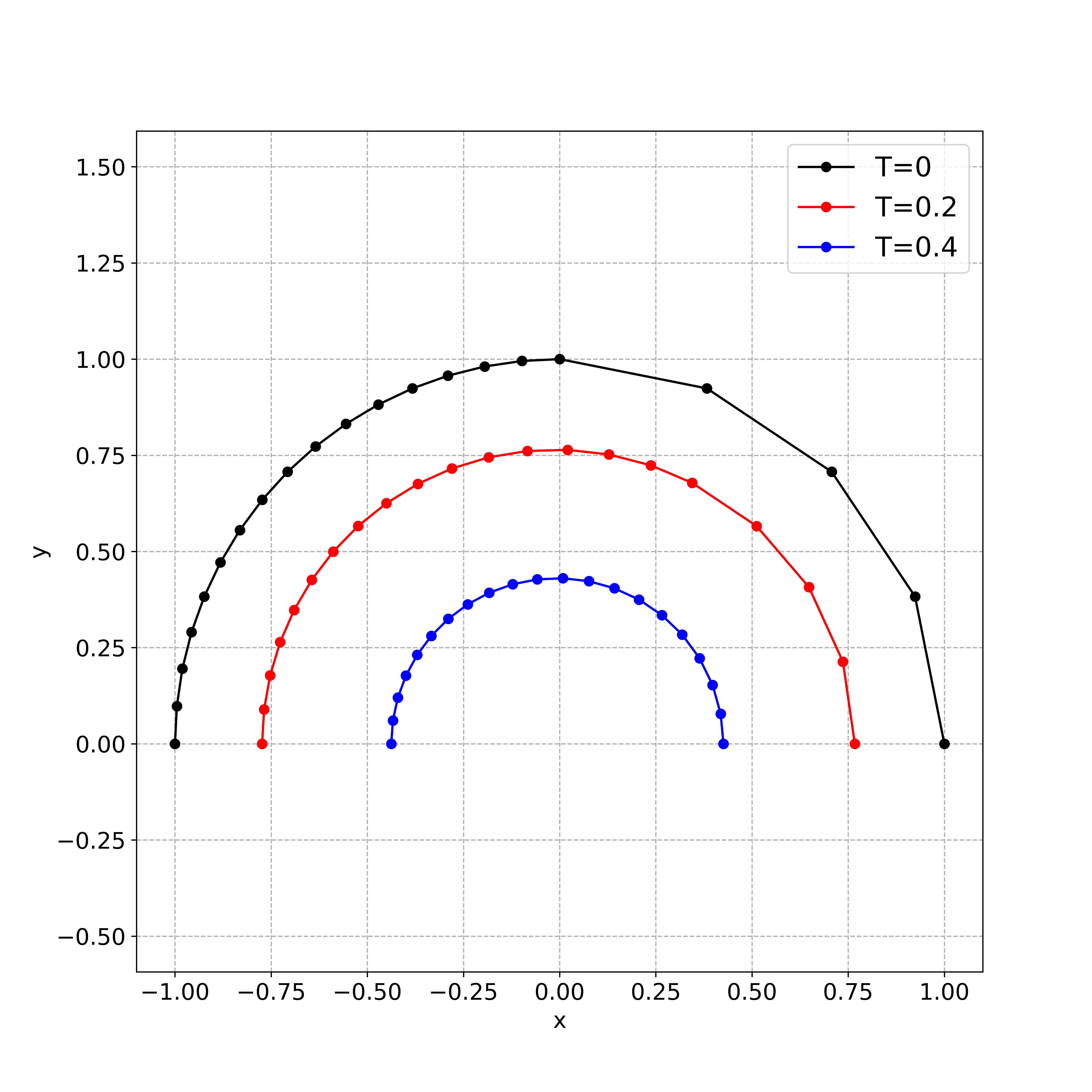}
            \caption{Curve evolution with nonuniform mesh (initial mesh size is small on the left and big on the right).}
            \label{Example_circle_que6_evol_1}
        \end{subfigure}
        \hfill
        \begin{subfigure}[b]{0.45\textwidth}
            \includegraphics[width=\textwidth]{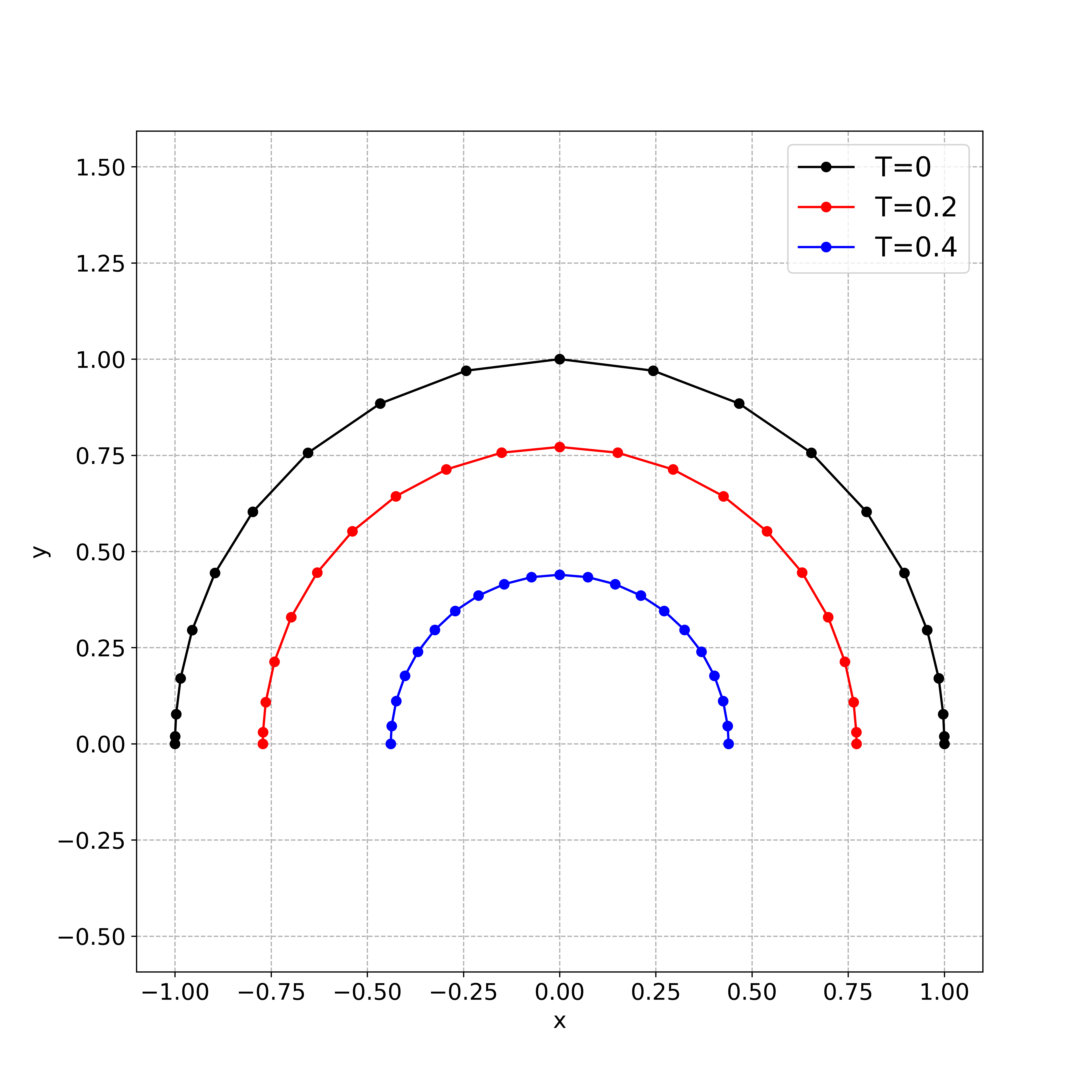}
           \caption{Curve evolution with nonuniform mesh (initial mesh size is small near the endpoints and big in the middle).}
            \label{Example_circle_temp_que6_evol}
        \end{subfigure}
        \caption{{Curve evolution for two different nonuniform meshes in Example~\ref{example_half_circle}.}}
        \label{example_circle_rates_que6_evol}
    \end{figure}

    \end{example}

    \end{example}
    
    \subsection{Model and governing equations in three dimensions (3D)}
    
    Following the framework in \cite{Bao2021} and \cite{Bao2023}, we consider a thin film on a flat substrate whose film–vapor interface, \(\Gamma(t)\), is parameterized by
    \[
    X(\cdot,t) = \bigl(X_1(\cdot,t),\, X_2(\cdot,t),\, X_3(\cdot,t)\bigr): \Gamma^0 \times [0,T] \rightarrow \mathbb{R}^3,
    \]
    with \(\Gamma^0\) denoting the initial configuration, and its evolution being described by the velocity field \(v(\cdot,t) = \partial_t X(\cdot,t)\). The film-substrate interface \(S_1(t)\) is a flat, two-dimensional domain that intersects \(\Gamma(t)\) along the contact line \(\partial\Gamma(t)=\Gamma(t)\cap S_1(t)\), which is assumed to form a simple closed curve oriented by the mapping \(X(p,t)\) for \(p\in\partial\Gamma^0\).
    
    The evolution of \(\Gamma(t)\) driven by surface diffusion, wherein the interface meets the substrate plane \(\mathbb{R}^2 \times \{0\}\) along \(\partial\Gamma(t)\), is governed by
    \begin{align*}
    v \cdot n &= \Delta_{\Gamma(t)} H \quad \text{on } \Gamma(t), \\
    H &= -\Delta_{\Gamma(t)} \mathrm{id} \cdot n \quad \text{on } \Gamma(t),
    \end{align*}
    where \(n\) denotes the unit normal vector and \(H\) is the mean curvature.
    
    For mean curvature flow, the evolution equations will be given as follows:
    \begin{align*}
    v \cdot n &= -H \quad \text{on } \Gamma(t), \\
    H &= -\Delta_{\Gamma(t)} \mathrm{id} \cdot n \quad \text{on } \Gamma(t).
    \end{align*}
    These equations are subject to the contact line condition
    \[
    X_3(\cdot,t)\big|_{\partial\Gamma} = X(\cdot,t)\cdot e_3 = 0 \quad \text{for } t\ge0,
    \]
    and the Young's law (contact angle condition)
    \[
    \mu_\partial \cdot n_\partial = \cos\theta.
    \]
    Here, \(\mu_\partial\) is the conormal vector on \(\partial\Gamma(t)\) (orthogonal to both \(\partial\Gamma(t)\) and \(n\)), \(n_\partial\) is the normal to \(\partial\Gamma(t)\) within the substrate plane, and \(\theta\) is the constant contact angle determined by the material properties of \(\Gamma(t)\) and the substrate. Moreover, for surface diffusion, the zero-mass flux condition
    \[
    \bigl(\mu_\partial \cdot \nabla_{\Gamma(t)} H\bigr)\big|_{\partial\Gamma} = 0 \quad \text{for } t\ge0,
    \]
    must also be imposed.

    For any test function 
\(
\eta := (\eta_1, \eta_2,\eta_3) \in H^1(\Gamma(t)) \times H^1(\Gamma(t)) \times  H^1(\Gamma(t)) ,
\)
the following equation holds obtained from integration by parts
\begin{align*}
    \int_{\Gamma(t)} Hn\cdot \eta &= \int_{\Gamma(t)} -\Delta_{\Gamma(t)} {\rm id} \cdot \eta  =  \int_{\Gamma(t)} \nabla_{\Gamma(t)} {\rm id} \cdot \nabla_{\Gamma(t)}\eta - \int_{\partial \Gamma(t)} \mu_{\partial}\cdot \eta \\
    & =  \int_{\Gamma(t)} \nabla_{\Gamma(t)} {\rm id} \cdot \nabla_{\Gamma(t)}\eta - \int_{\partial \Gamma(t)} (\mu_{\partial}\cdot n_\partial) ( \eta  \cdot n_\partial ) - \int_{\partial \Gamma(t)} (\mu_{\partial} \cdot e_3) (\eta\cdot e_3)  \\
    & = \int_{\Gamma(t)} \nabla_{\Gamma(t)} {\rm id} \cdot \nabla_{\Gamma(t)}\eta - \int_{\partial \Gamma(t)} \cos\theta \,n_\partial \cdot \eta - \int_{\partial \Gamma(t)} \sin\theta \,e_3 \cdot \eta  ,  
\end{align*}
where the last equality is obtained by applying Young's law.

Then the functional \( L_h^{m-1}(\eta_h) \) is defined on the finite element function space \((S_h^{m-1})^3\) as follows:
\[
L_h^{m-1}(\eta_h) :=\int_{\Gamma_h^{m-1}} \nabla_{\Gamma_h^{m-1}} {\rm id} \cdot \nabla_{\Gamma_h^{m-1}} \eta_h - \int_{\partial \Gamma_h^{m-1}}^{(h)} \cos \theta \, {{n_{\partial,h}^{m-1}}}\cdot \eta_h -\int_{\partial \Gamma_h^{m-1}}^{(h)} \sin \theta \, e_3\cdot \eta_h,
\]
{where $n_{\partial,h}^{m-1}$ denotes the piecewise normal vector on the discrete curve  $\partial \Gamma_h^{m-1}$ within the substrate plane.}
By virtue of the Riesz representation theorem, it is deduced that there exists \(\nu_h^{m-1} \in (S_h^{m-1} )^3\) such that
\[
L_h^{m-1}(\eta_h) = \bigl( \nu_h^{m-1}, \eta_h \bigr)^{(h)} \quad \text{for all } \eta_h \in (S_h^{m-1})^3.
\]

Then the tangential vector \( T_h^{m-1} \) is defined similarly to the equation \eqref{eq:T_h^m-higher-order-2D}, and it is required to determine \((T_h^{m-1}, \lambda_h) \in (S_h^{m-1})^3 \times S_h^{m-1}\) such that
\begin{equation}\label{eq:T_h^m-higher-order-3D-1}
\begin{aligned}
\bigl( T_h^{m-1}, \eta_h \bigr)^{(h)} + \bigl( \lambda_h\, n_h^{m-1}, \eta_h \bigr)^{(h)} &= \bigl( \nu_h^{m-1}, \eta_h \bigr)^{(h)}, \quad \forall \eta_h \in (S_h^{m-1})^3, \\[1mm]
\bigl( T_h^{m-1}, n_h^{m-1}\, \chi_h \bigr)^{(h)} &= 0, \quad \forall \chi_h \in S_h^{m-1}.
\end{aligned}
\end{equation}

    \subsection{Numerical scheme for mean curvature flow and surface diffusion in 3D}

    Let \( \mathring{S}_h(\Gamma_h^{m-1}) \) denote the subspace of \( S_h^{m-1} \) consisting of finite element functions that vanish on \( \partial \Gamma_h^{m-1} \). The product space \( \mathbf{X}_h^{m-1} := S_h^{m-1} \times S_h^{m-1} \times \mathring{S}_h(\Gamma_h^{m-1}) \) is then defined. The discrete velocity \( v_h^m \) of the flow map \( X_h^m \) is assumed to lie within \( \mathbf{X}_h^{m-1} \), ensuring that the contact line condition is satisfied.

Using the vector \( T_h^{m-1} \) defined in \eqref{eq:T_h^m-higher-order-3D-1}, the BGN-MDR scheme for mean curvature flow of an open surface with a moving contact line is proposed. The goal is to find \( (v_h^m, \lambda_h^m, c^m) \in \mathbf{X}_h^{m-1} \times S_h^{m-1} \times \mathbb{R} \) such that the following system is satisfied:
\begin{subequations}\label{eq:MCF-num-higher-3D}
    \begin{align}
        & \int_{\Gamma_h^{m-1}} \nabla_{\Gamma_h^{m-1}} (\tau v_h^m + {\rm id}) \cdot \nabla_{\Gamma_h^{m-1}} \eta_h - (  \lambda_h^m {n}_h^{m-1}, \eta_h)^{(h)} - c^m (T_h^{m-1}, \eta_h)^{(h)} \nonumber \\
        & \qquad = \int_{\partial \Gamma_h^{m-1}} \cos(\theta) \, {n}_{\partial,h}^{m-\frac{1}{2}} \cdot \eta_h, \quad \forall \eta_h \in \mathbf{X}_h^{m-1}, \label{eq:MCF-num-1-higher-3D} \\
        & (v_h^m \cdot {n}_h^{m-1}, \chi_h)^{(h)} + (\lambda_h^m, \chi_h)^{(h)} = 0, \quad \forall \chi_h \in S_h^{m-1}, \label{eq:MCF-num-2-higher-3D} \\
        & (v_h^m ,T_h^{m-1})^{(h)} + \alpha c^m \|T_h^{m-1}\|_{L^2_h} = 0. \label{eq:MCF-num-3-higher-3D}
    \end{align}
\end{subequations}
For the numerical approximation of \( n_\partial \) on the discrete cure \( \partial \Gamma_h^{m-1} \), we adopt the approach presented in equation (3.7) of \cite{Bao2023}. This choice is not arbitrary; it is selected to ensure the dissipation of surface energy. The approximated conormal vector \( {n}_\partial^{m-\frac{1}{2}} \) at the boundary \( \partial \Gamma_h^{m-1} \) is defined by
    \begin{align}\label{n-partial-def}
        {n}_\partial^{m-\frac{1}{2}} = \frac{1}{2} \left( \partial_s X_h^{m-1}  + \partial_s X_h^m  \right) \times e_3,
    \end{align}
    where \( s \) denotes the arc length parameter on the boundary curve \( \partial \Gamma_h^{m-1} \).

Additionally, the BGN-MDR numerical scheme for surface diffusion of an open surface with a moving contact line is formulated as follows: find \( (v_h^m, \lambda_h^m,  c^m) \in \mathbf{X}_h^{m-1} \times S_h^{m-1} \times \mathbb{R} \) such that the system
\begin{subequations}\label{eq:SD-num-higher-3D}
    \begin{align}
        & \int_{\Gamma_h^{m-1}} \nabla_{\Gamma_h^{m-1}} (\tau v_h^m + {\rm id}) \cdot \nabla_{\Gamma_h^{m-1}} \eta_h - (  \lambda_h^m {n}_h^{m-1}, \eta_h)^{(h)} - c^m (T_h^{m-1} , \eta_h)^{(h)} \nonumber \\
        & \qquad = \int_{\partial \Gamma_h^{m-1}} \cos(\theta) \, {n}_{\partial,h}^{m-\frac{1}{2}} \cdot \eta_h, \quad \forall \eta_h \in \mathbf{X}_h^{m-1}, \label{eq:SD-num-1-higher-3D} \\
        & (v_h^m \cdot {n}_h^{m-1}, \chi_h)^{(h)} + \int_{\Gamma_h^{m-1}} \nabla_{\Gamma_h^{m-1}} \lambda_h^m \cdot \nabla_{\Gamma_h^{m-1}} \chi_h = 0, \quad \forall \chi_h \in S_h^{m-1}, \label{eq:SD-num-2-higher-3D} \\
        & (v_h^m , T_h^{m-1})^{(h)} + \alpha c^m \|T_h^{m-1}\|_{L^2_h} = 0 \label{eq:SD-num-3-higher-3D}
    \end{align}
\end{subequations}
is satisfied.


{
\begin{remark}\upshape
\label{remark-well-posedness-3D}
In three dimensions, the appearance of $X_h^m$ in the definition of $n_\partial^{m-\frac{1}{2}}$ in \eqref{n-partial-def}, which is subsequently used in \eqref{eq:MCF-num-1-higher-3D} and \eqref{eq:SD-num-1-higher-3D}, prevents the establishment of well-posedness as in the two-dimensional case. This is due to the contribution of the boundary term in the bilinear form. While this boundary term is crucial for preserving energy dissipation, it precludes a formal proof of well-posedness. However, if $n_\partial^{m-1}$ is used in place of $n_\partial^{m-\frac{1}{2}}$ in \eqref{eq:MCF-num-1-higher-3D} and \eqref{eq:SD-num-1-higher-3D}, the boundary term can instead be treated as part of the linear form. This modification enables the establishment of well-posedness by employing an argument analogous to the proof of well-posedness in the two-dimensional case. 
\end{remark}
}

Furthermore, the BGN-MDR numerical scheme for open surfaces with moving contact lines preserves unconditional energy stability for both mean curvature flow and surface diffusion due to the specific construction of \( {n}_{\partial,h}^{m-\frac{1}{2}} \), as outlined in \cite{Bao2023}.

\begin{theorem}\label{thm:area-decreasing-mass-sd-3D}
    Let \( (v_h^m, \lambda_h^m, c^m) \in \mathbf{X}_h^{m-1} \times S_h^{m-1} \times \mathbb{R} \) be the solution to the weak formulation \eqref{eq:SD-num-higher-3D}. The scheme guarantees that the discrete surface energy \( W_h^m := |\Gamma_h^m| - \cos \theta |S_1^m| \) is non-increasing over time, i.e.,
    \begin{align}\label{thm:area-decreasing-mass-sd-1-3D}
        W_h^m \le W_h^{m-1} \quad \text{for } m = 1, 2, \ldots, N.
    \end{align}
Here, \(S_1^m\) denotes the discrete film-substrate interface at the time level \(t^m\), which serves as an approximation to the film-substrate interface \(S_1(t^m)\).
\end{theorem}

\begin{proof}
By choosing \( \eta_h := v_h^m\) in \eqref{eq:SD-num-1-higher-3D} and \( \chi_h := \lambda_h^m \) in \eqref{eq:SD-num-2-higher-3D}, and sum up equations \eqref{eq:SD-num-1-higher-3D}--\eqref{eq:SD-num-3-higher-3D}, we obtain the following relation:
    \begin{equation}\label{thm:area-decreasing-mass-SD-2-3D}
        \begin{aligned}
            &\int_{\Gamma_h^{m-1}} \nabla_{\Gamma_h^{m-1}} (\tau v_h^m + {\rm id}) \cdot \nabla_{\Gamma_h^{m-1}}v_h^m + \|\nabla_{\Gamma_h^{m-1}} \lambda_h^m\|_{L^2(\Gamma_h^{m-1})}^2 + \frac{1}{\alpha}s(c^m)^2 \|T_h^m\|_{L^2_h} \\
            &= \frac{1}{\tau} \cos \theta \int_{\partial \Gamma_h^{m-1}} n_{\partial,h}^{m-\frac{1}{2}} \cdot (X_h^m - {\rm id}).
        \end{aligned}
    \end{equation}
From \cite[Lemma 3.1]{Bao2023} we know that the right-hand side of \eqref{thm:area-decreasing-mass-SD-2-3D} can be written as
    \begin{align}\label{bottom-area-3D}
        \int_{\partial \Gamma_h^{m-1}} {n}_{\partial,h}^{m-\frac{1}{2}} \cdot (X_h^m - {\rm id}) = |S_1^m| - |S_1^{m-1}| . 
    \end{align}
Therefore, substituting \eqref{bottom-area-3D} and \eqref{area-decreasing-important-property} into \eqref{thm:area-decreasing-mass-SD-2-3D} leads to $W_h^m - W_h^{m-1} \le 0 $. 
\hfill\end{proof}

    The convergence and performance of the BGN-MDR method for mean curvature flow of open surfaces with moving contact lines are shown in the following two examples. 

\begin{example}\label{example_half_sphere}\upshape
We simulate the mean curvature flow of half-sphere with a moving contact line on the plane $z=0$ by using the BGN-MDR  scheme in \eqref{eq:MCF-num-higher-3D}. The exact solution of this problem is a half-sphere with radius $r(t)=\sqrt{1-4t}$.

Figures~\ref{example_sphere_rates}(a) and \ref{example_sphere_rates}(b) show that the errors from spatial and temporal discretizations are \(O(h^2)\) and \(O(\tau)\), respectively, in the \(L^\infty(0,0.1;L^2)\)-norm, demonstrating the convergence of the proposed BGN-MDR scheme for open surfaces with moving contact lines. 

In this example, we also compare the results obtained with different values of $\alpha$. Figure~\ref{Example_sphere_curvature} presents the evolution of discrete curvature for various choices of $\alpha$, while Figure~\ref{Example_sphere_energy} depicts the corresponding evolution of discrete energy. It can be observed that the curves for the three values of $\alpha$ nearly overlap, indicating that discrete curvature and energy can be computed with good accuracy for a wide range of $\alpha$. This demonstrates that $\alpha$ can be chosen freely to optimize the mesh quality.

    \begin{figure}[htbp]
        \centering
    
        \begin{subfigure}[b]{0.45\textwidth}
            \includegraphics[width=\textwidth]{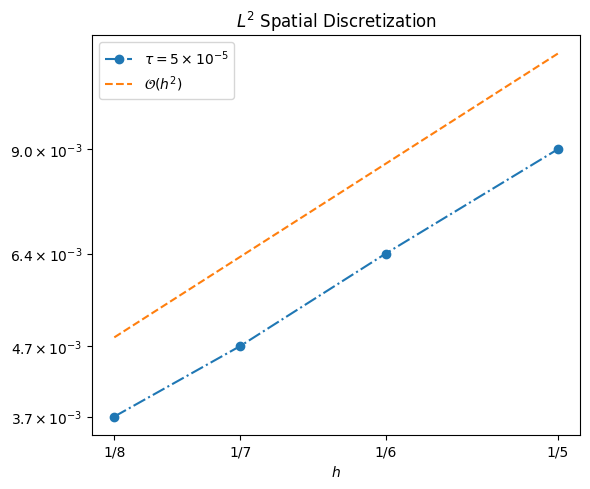}
            \caption{Spatial discretization errors}
            \label{Example_sphere_spat}
        \end{subfigure}
        \hfill
        \begin{subfigure}[b]{0.45\textwidth}
            \includegraphics[width=\textwidth]{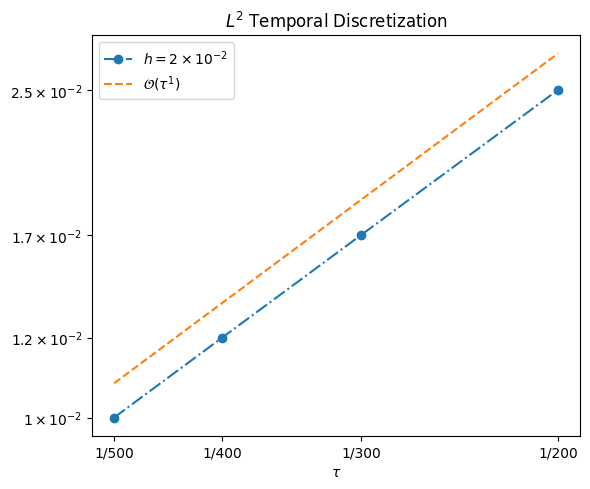}
            \caption{Temporal discretization errors}
            \label{Example_sphere_temp}
        \end{subfigure}
        \caption{\(L^\infty(0,0.1;L^2)\) errors of the numerical solutions in Example~\ref{example_half_sphere}}
        \label{example_sphere_rates}
    \end{figure}

    \begin{figure}[htbp]
        \centering
    
        \begin{subfigure}[b]{0.45\textwidth}
            \includegraphics[width=\textwidth]{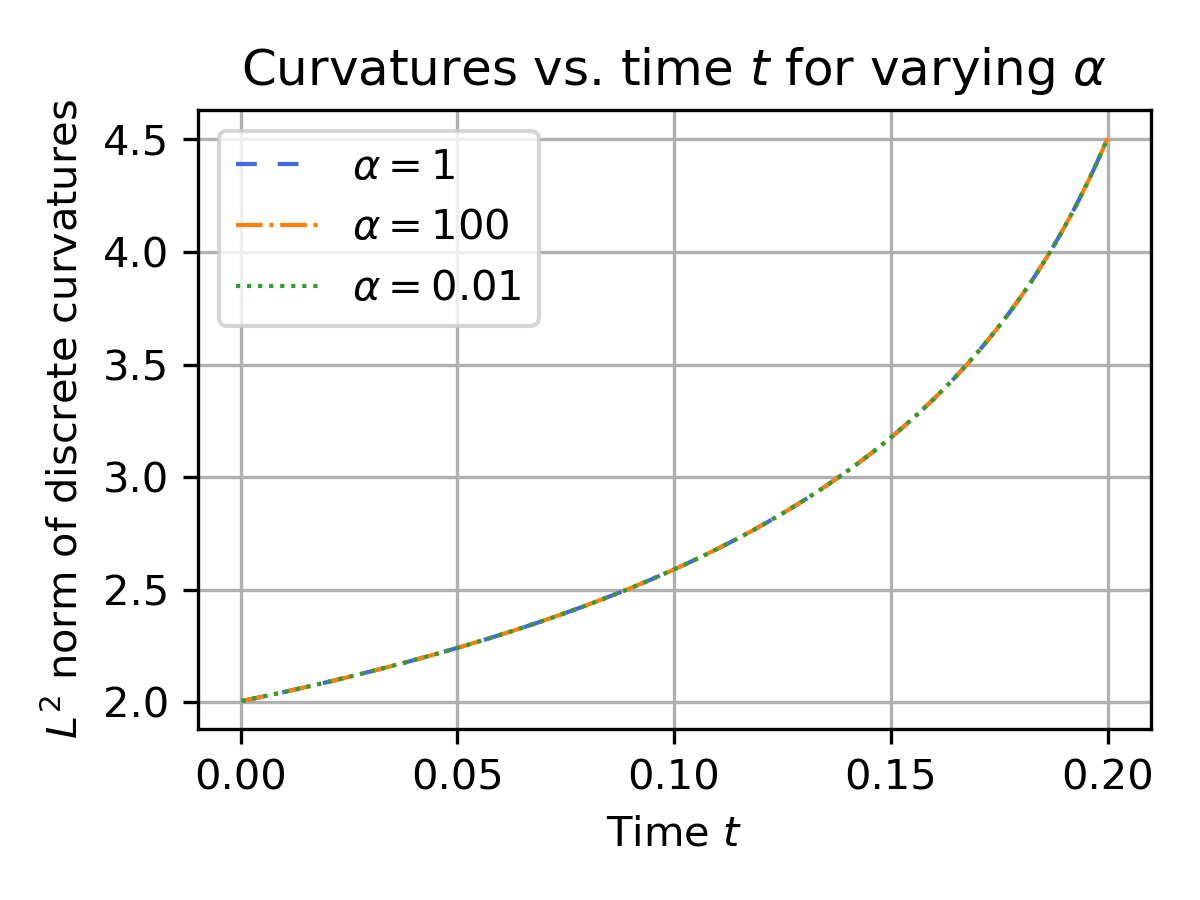}
            \caption{Discrete curvature for varying $\alpha$}
            \label{Example_sphere_curvature}
        \end{subfigure}
        \hfill
        \begin{subfigure}[b]{0.45\textwidth}
            \includegraphics[width=\textwidth]{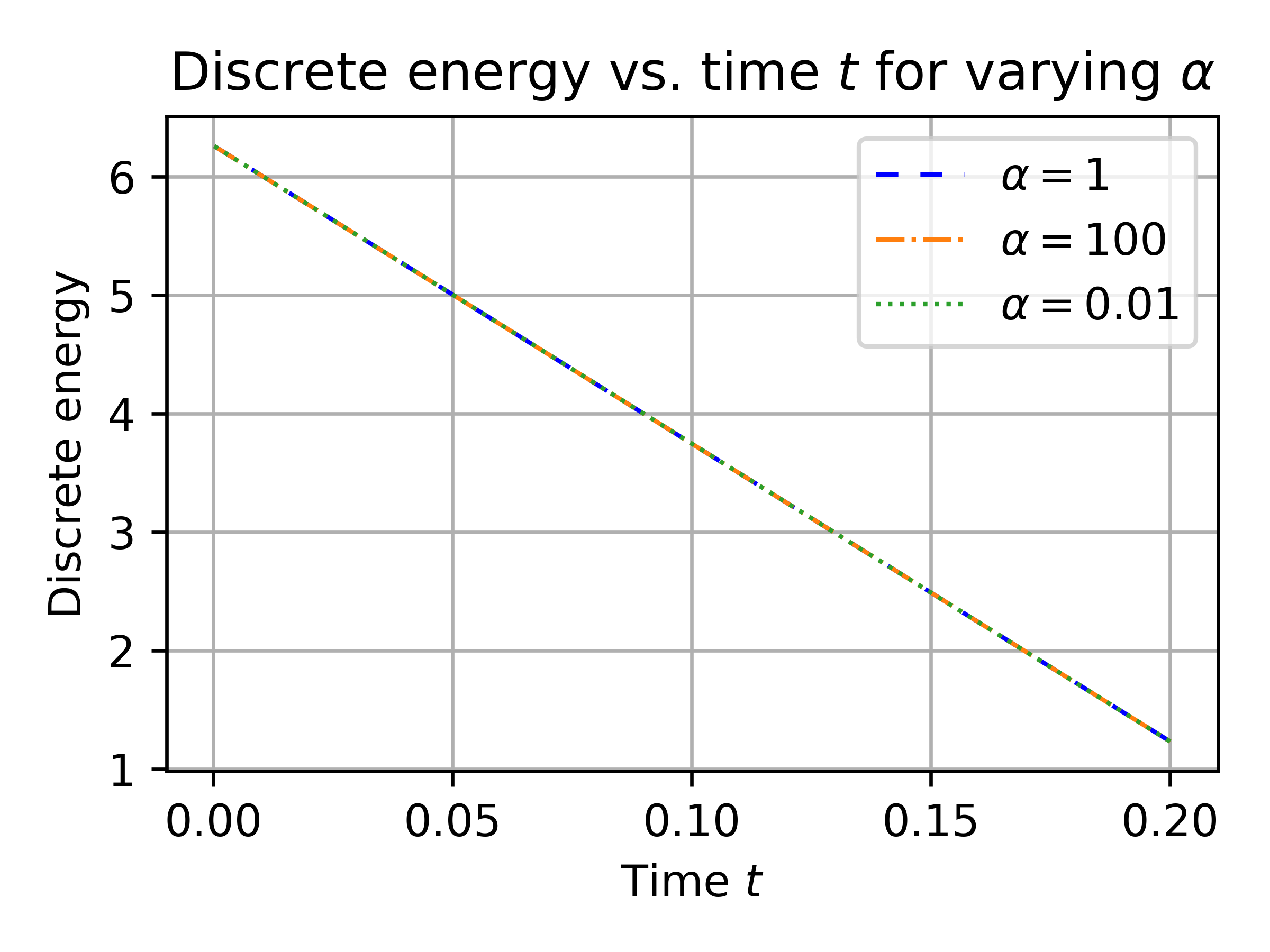}
            \caption{Discrete energy for varying $\alpha$}
            \label{Example_sphere_energy}
        \end{subfigure}
        \caption{{Discrete energy and curvature versus time for varying $\alpha$ in Example~\ref{example_half_sphere}.}}
        \label{example_sphere_ref}
    \end{figure}
    
    \end{example}

\begin{example}\label{Example7}\upshape
Next, we present numerical simulations for surface diffusion of an open surface with a moving contact line constrained to the plane $z=0$ with $60^\circ$ and $120^\circ$ contact angles, respectively, with the initial surface being a $1 \times 6 \times 1$ box centered at $(0,0,0)$. We present numerical simulations with the BGN method and the BGN-MDR method with parameter $\alpha = 0.01$.

The numerical simulation with the $60^\circ$ contact angle is presented in Figure \ref{eg:Example7_60}. By using mesh size $h = 0.2$ and time step size $\tau = 10^{-2}$, both the BGN and BGN-MDR schemes exhibit satisfactory mesh quality, as illustrated in Figures~\ref{Example7d_60} and~\ref{Example7e_60}. In contrast, when a smaller time step size $\tau = 10^{-3}$ is used, the BGN scheme becomes unstable and breaks down at $T = 0.549$, and the BGN-MDR scheme remains stable and yields good mesh quality throughout the evolution; see Figures~\ref{Example7b_60} and~\ref{Example7c_60}.

\begin{figure}[htbp]
            \centering
 %
            \begin{subfigure}[b]{0.5\textwidth}
            \centering
                \includegraphics[width=0.85\textwidth]{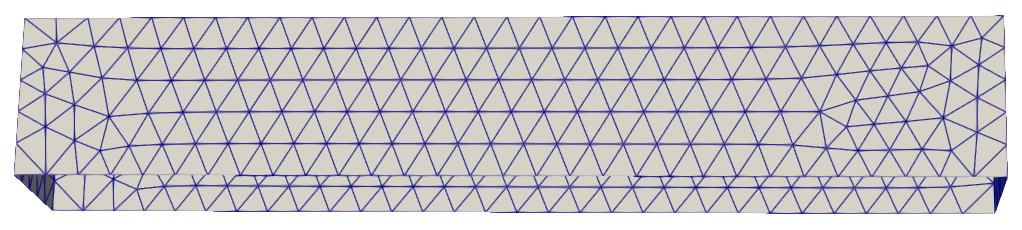}
                \caption{Initial Surface}
                \label{Example7a_60}
            \end{subfigure}
        
            \vspace{10pt}
        
            \begin{subfigure}[b]{0.45\textwidth}
            \centering
                \includegraphics[width=0.75\textwidth]{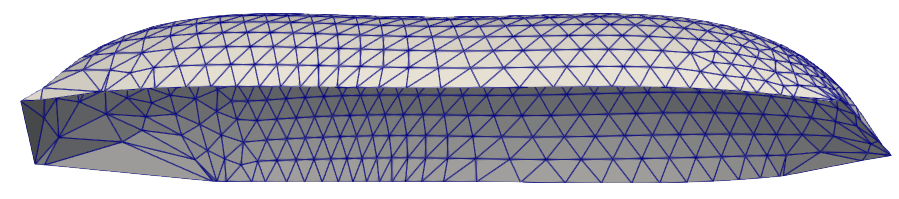}
                \caption{BGN with $\tau=10^{-3}$ at $T = 0.454$}
                \label{Example7b_60}
            \end{subfigure}
            \hspace{10pt}
            \begin{subfigure}[b]{0.45\textwidth}
            \centering
                \includegraphics[width=0.6\textwidth]{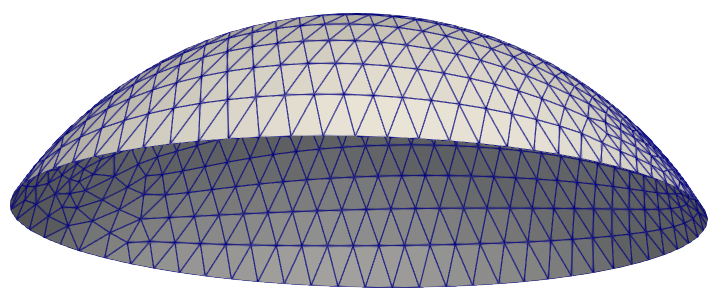}
                \caption{BGN-MDR with $\tau=10^{-3}$ at $T=2$}
                \label{Example7c_60}
            \end{subfigure}
        
            \vspace{10pt}
        
            \begin{subfigure}[b]{0.45\textwidth}
            \centering
                \includegraphics[width=0.6\textwidth]{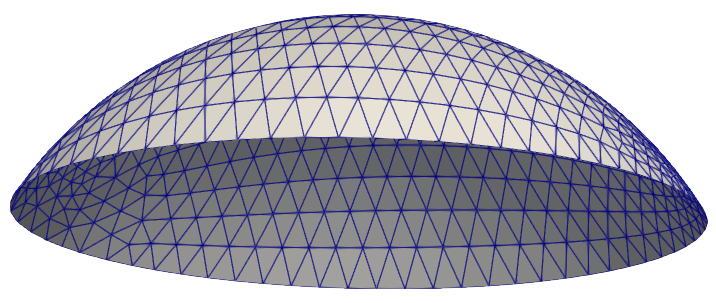}
                \caption{BGN with $\tau=10^{-2}$ at $T=2$}
                \label{Example7d_60}
            \end{subfigure}
            \hspace{10pt}
            \begin{subfigure}[b]{0.45\textwidth}
            \centering
                \includegraphics[width=0.6\textwidth]{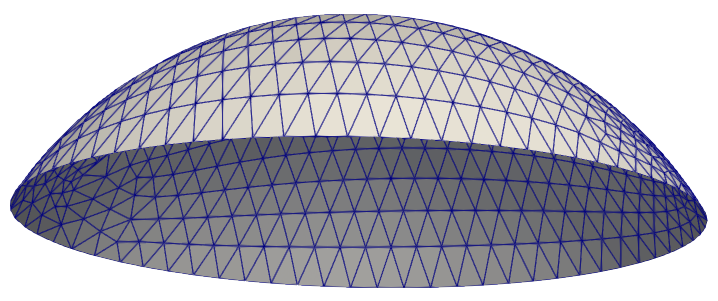}
                \caption{BGN-MDR with $\tau=10^{-2}$ at $T=2$}
                \label{Example7e_60}
            \end{subfigure}
        
            \vspace{-5pt}
            \caption{Surface evolution with a $60^\circ$ contact angle in Example~\ref{Example7}
            }
            \label{eg:Example7_60}
        \end{figure}
        
The numerical simulation with the $120^\circ$ contact angle in presented in Figure \ref{eg:Example7_120}. By using mesh size $h = 0.2$ and time step size $\tau = 10^{-2}$, both schemes maintain satisfactory mesh quality, as shown in Figures~\ref{Example7d_120} and~\ref{Example7e_120}. However, when a smaller time step size $\tau = 10^{-3}$ is used, the BGN scheme becomes unstable and breaks down at $T = 0.227$. In contrast, the BGN-MDR scheme remains stable and continues to produce good mesh quality throughout the evolution up to $T = 2$; see Figures~\ref{Example7b_120} and~\ref{Example7c_120}.

        \begin{figure}[htbp]
            \centering
 %
%
        
            \begin{subfigure}[b]{0.45\textwidth}
            \centering
                \includegraphics[width=0.75\textwidth]{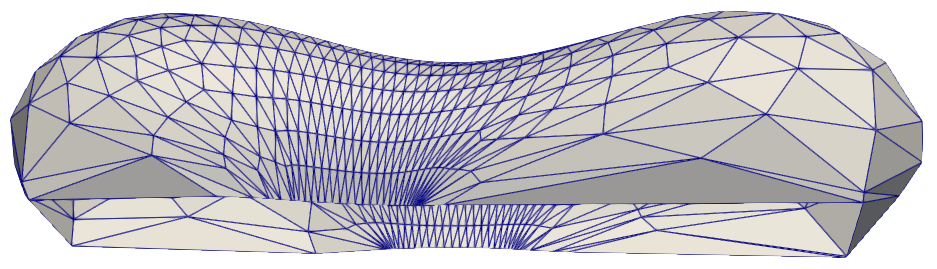}
                \caption{BGN with $\tau=10^{-3}$ at $T=0.181$}
                \label{Example7b_120}
            \end{subfigure}
            \hspace{10pt}
            \begin{subfigure}[b]{0.45\textwidth}
            \centering
                \includegraphics[width=0.5\textwidth]{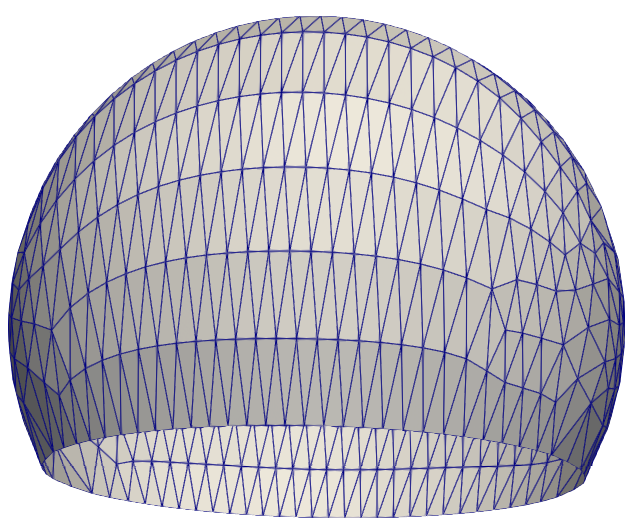}
                \caption{BGN-MDR with $\tau=10^{-3}$ at $T=2$}
                \label{Example7c_120}
            \end{subfigure}
        
            \vspace{10pt}
        
            \begin{subfigure}[b]{0.45\textwidth}
            \centering
                \includegraphics[width=0.5\textwidth]{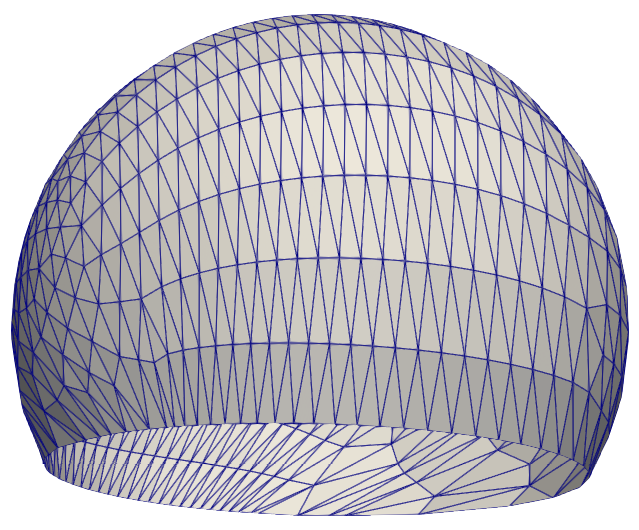}
                \caption{BGN with $\tau=10^{-2}$ at $T=2$}
                \label{Example7d_120}
            \end{subfigure}
            \hspace{10pt}
            \begin{subfigure}[b]{0.45\textwidth}
            \centering
                \includegraphics[width=0.5\textwidth]{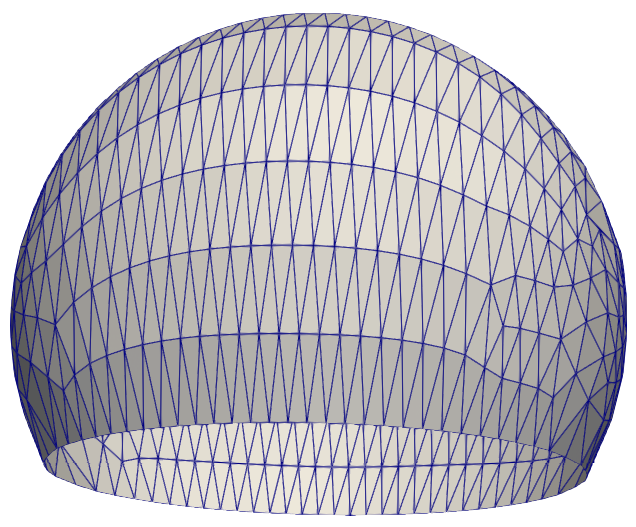}
                \caption{BGN-MDR with $\tau=10^{-2}$ at $T=2$}
                \label{Example7e_120}
            \end{subfigure}
        
            \vspace{-5pt}
            \caption{Surface evolution with a $120^\circ$ contact angle in Example~\ref{Example7}
            }
            \label{eg:Example7_120}
        \end{figure}

    \end{example}

\section{Conclusion}
We propose a new parametric FEM, referred to as the BGN-MDR method, for simulating the evolution of curves and surfaces under mean curvature flow and surface diffusion, applicable to both closed surfaces and open surfaces with moving contact lines or points. The proposed BGN-MDR method bridges the strengths of the BGN and MDR methods: it maintains good mesh quality similar to the MDR method, while also preserving energy stability as in the BGN method. Numerical simulations demonstrate the advantages of the proposed approach on benchmark problems involving both closed and open surfaces.

\renewcommand{\refname}{\bf References}
\bibliographystyle{abbrv}
\bibliography{main}

\end{document}